\documentclass[aap,preprint]{imsart}

\usepackage{enumerate}
\usepackage{amssymb}
\usepackage{amsthm,amsmath}
\usepackage{amsfonts}
\usepackage{mathrsfs}
\usepackage{graphicx}
\usepackage{mathrsfs}
\usepackage{stmaryrd}
\usepackage{bbm}
\usepackage{cases}
\usepackage{amsfonts}
\usepackage{txfonts}

\arxiv{ arXiv:1406.7446}

\startlocaldefs

\newcommand{\be}{\begin{eqnarray}}
\newcommand{\ee}{\end{eqnarray}}
\newcommand{\ce}{\begin{eqnarray*}}
\newcommand{\de}{\end{eqnarray*}}
\newtheorem{theorem}{Theorem}[section]
\newtheorem{lemma}[theorem]{Lemma}
\newtheorem{remark}[theorem]{Remark}
\newtheorem{definition}[theorem]{Definition}
\newtheorem{proposition}[theorem]{Proposition}
\newtheorem{Examples}[theorem]{Example}
\newtheorem{corollary}[theorem]{Corollary}

\def\v{{\mathbf{v}}}
\def\e{{\mathrm{e}}}
\def\eps{\varepsilon}

\def\p{\partial}

\def\[{{\Big[}}
\def\]{{\Big]}}
\def\<{{\langle}}
\def\>{{\rangle}}
\def\({{\Big(}}
\def\){{\Big)}}

\def\bx{{\mathbf{x}}}

\def\dif{{\mathord{{\rm d}}}}

\def\no{\nonumber}
\def\={&\!\!=\!\!&}
\def\bt{\begin{theorem}}
\def\et{\end{theorem}}
\def\bl{\begin{lemma}}
\def\el{\end{lemma}}
\def\br{\begin{remark}}
\def\er{\end{remark}}

\def\bd{\begin{definition}}
\def\ed{\end{definition}}
\def\bp{\begin{proposition}}
\def\ep{\end{proposition}}
\def\bc{\begin{corollary}}
\def\ec{\end{corollary}}
\def\bx{\begin{Examples}}
\def\ex{\end{Examples}}

\def\cM{{\mathcal M}}

\def\cT{{\mathcal T}}

\def\mE{{\mathbb E}}

\def\mH{{\mathbb H}}
\def\mI{{\mathbb I}}

\def\mL{{\mathbb L}}
\def\mM{{\mathbb M}}
\def\mN{{\mathbb N}}

\def\mR{{\mathbb R}}

\def\mT{{\mathbb T}}

\def\mW{{\mathbb W}}

\def\bP{{\mathbf P}}

\def\sC{{\mathscr C}}

\def\sF{{\mathscr F}}

\def\sV{{\mathscr V}}

\def\geq{\geqslant}
\def\leq{\leqslant}

\def\div{\mathord{{\rm div}}}

\def\bP{{\mathbf P}}

\def\bK{{\mathbf K}}

\def\u{{\mathbf u}}

\endlocaldefs

\begin{document}

\begin{frontmatter}

\title{Stochastic differential equations  with Sobolev diffusion and singular drift}
\runtitle{SDEs with Sobolev diffusion and singular drift}

\author{\fnms{Xicheng} \snm{Zhang}\ead[label=e2]{XichengZhang@gmail.com}}
\address{School of Mathematics and Statistics,\\ Wuhan University,\\ Wuhan, Hubei 430072, P.R.China\\ \printead{e2}}
\affiliation{Wuhan University}

\runauthor{}

\begin{abstract}
In this paper we study properties of solutions to stochastic differential equations
with Sobolev diffusion coefficients and singular drifts. The properties we study include
stability with respect to the coefficients, weak differentiability with respect to starting points,
and the Malliavin differentiability with respect to sample paths.
We also establish Bismut-Elworthy-Li's formula for the solutions.
As an application, we use the stochastic Lagrangian representation of
incompressible Navier-Stokes equations
given by Constantin-Iyer \cite{Co-Iy} to prove the local well-posedness of
NSEs in $\mR^d$ with initial values in the first order Sobolev space
$\mW^1_p(\mR^d;\mR^d)$ provided $p>d$.
\end{abstract}



\end{frontmatter}

\section{Introduction and Main Results}

Consider the following stochastic differential equation (abbreviated as SDE) in $\mR^d$:
\begin{align}
\dif X_t=b_t(X_t)\dif t+\dif W_t, \ \ t\geq 0,\ \ X_0=x\in\mR^d,\label{KL1}
\end{align}
where $(W_t)_{t\geq 0}$ is a $d$-dimensional standard Brownian motion on some probability space $(\Omega,\sF,P)$.
It is a classical result due to Veretennikov \cite{Ve} that when $b$ is bounded and Borel measurable,
the SDE above admits a unique strong solution. Furthermore,
for almost all $\omega$, the following random ordinary differential equation
$$
\dif X_t(\omega)=b_t(X_t(\omega)+W_t(\omega))\dif t, \ \ t\geq 0,\ \ X_0=x
$$
has a unique solution (cf. Davie \cite{Da}).
Recently, in \cite{Me} and \cite{Mo-Ni-Pr}, the Malliavin and Sobolev differentiabilities of $X_t(x,\omega)$
with respect to the sample path $\omega$
and with respect to  the starting point $x$ were studied,
and these differentiabilities were used  to study stochastic transport equations.
In a remarkable paper \cite{Kr-Ro}, Krylov and R\"ockner proved the existence and
uniqueness of strong solutions to SDE (\ref{KL1}) under the assumption
$$
b\in L^q(\mR_+; L^p(\mR^d)\ \mbox{ with $p,q\in(1,\infty)$ and }\tfrac{d}{p}+\tfrac{2}{q}<1,
$$
by using the Girsanov transformation and some estimates from the theory of PDEs.
Subsequently, the results of \cite{Kr-Ro} were extended to the case of multiplicative noises in
\cite{Zh2} (see also \cite{Gy-Ma, Zh0} for related results).
The Sobolev differentiability of solutions was also obtained in \cite{Fe-Fl0, Fe-Fl}.
The recent interest in studying the Sobolev differentiability for (\ref{KL1}) with singular drift is
partly due to the discovery of Flandoli, Gubinelli and Priola \cite{Fl-Gu-Pr} that
noises can prevent the singularity for linear transport equations (see also \cite{Fe-Fl0}).

In this paper we consider the following SDE: for given $T<S$,
\begin{align}
\dif X_{t,s}=b_s(X_{t,s})\dif s+\sigma_s(X_{t,s})\dif W_s,\ X_{t,t}=x, \ T\leq t\leq s\leq S,\label{SDE00}
\end{align}
where $b:[T,S]\times\mR^d\to\mR^d$ and $\sigma:[T,S]\times\mR^d\to\mM^d$ are two Borel functions,
and $(W_s)_{s\in[T,S]}$
is a $d$-dimensional standard Brownian motion on the classical Wiener space $(\Omega,\sF,P;\mH)$.  Here, $\mM^d$ denotes the set of all $d\times d$-matrices,
$\Omega$ is the space of all continuous functions from $[T,S]$ to $\mR^d$, $\sF$ is the Borel-$\sigma$ field,  $P$ is the Wiener measure, and $\mH\subset\Omega$
is the Cameron-Martin space. We make the following assumption on $\sigma$:
\begin{enumerate}[{\bf (H$^\alpha_K$)}]
\item
there exist constants $K\geq 1$ and $\alpha\in (0, 1)$ such that for all $(t,x)\in[T,S]\times\mR^d$,
\begin{align}
K^{-1}|\xi|\leq |\sigma^{\textrm{t}}_t(x)\xi|\leq K|\xi|,\ \xi\in\mR^d,\label{Con2}
\end{align}
and for all $t\in[T,S]$ and $x,y\in\mR^d$,
$$
\|\sigma_t(x)-\sigma_t(y)\|\leq K|x-y|^\alpha.
$$
Here and in the remainder of this paper, $\sigma^{\mathrm{t}}$ denotes the transpose of matrix $\sigma$,
$|\cdot|$ the Euclidiean norm and $\|\cdot\|$ the Hilbert-Schmidt norm.
\end{enumerate}

Throughout this work, for simplicity of presentation, we assume $S-T\leq 1$ 
so that all the constants appearing below are independent of the length of the time interval $[T,S]$.
Our main result of this paper is:
\bt\label{Main}
Assume that $\sigma$ satisfies {\bf (H$^\alpha_K$)}. Suppose also that one of the following two conditions holds:
\begin{enumerate}[(i)]
\item $\sigma_t(x)=\sigma_t$ is independent of $x$ and for some $p,q\in(1,\infty)$ with $\tfrac{d}{p}+\tfrac{2}{q}<1$,
$$
b\in L^q([T,S]; L^p(\mR^d))=:\mL^q_p(T,S).
$$
\item $\nabla\sigma, b\in \mL^q_p(T,S)$ for some $q=p>d+2$.
\end{enumerate}
Then we have the following conclusions:
\begin{enumerate}[{\bf (A)}]
\item For any $(t,x)\in[T,S]\times\mR^d$, there is a unique strong solution denoted by $X_{t,s}(x)$ or $X^{b,\sigma}_{t,s}(x)$ to SDE (\ref{SDE00}),
which has a jointly continuous version with respect to $s$ and $x$.
\item For each $s\geq t$ and almost all $\omega$, $x\mapsto X_{t,s}(x,\omega)$ is weakly differentiable.
Furthermore, for any $p'\geq 1$, the Jacobian matrix $\nabla X_{t,s}(x)$ satisfies
\begin{align}\label{EU44}
\begin{split}
&\mathrm{ess.}\sup_{x\in\mR^d}\mE\left(\sup_{s\in[t,S]}|\nabla X_{t,s}(x)|^{p'}\right)\\
&\quad\leq C=C\Big(d,p,q,K,\alpha,p', \|b\|_{\mL^q_p(t,S)}, \|\nabla\sigma\|_{\mL^q_p(t,S)}\Big),
\end{split}
\end{align}
where the constant $C$ is increasing with respect to $\|b\|_{\mL^q_p(t,S)}$ and $\|\nabla\sigma\|_{\mL^q_p(t,S)}$.
\item For each $s\geq t$ and $x\in\mR^d$, the random variable $\omega\mapsto X_{t,s}(x,\omega)$ is Malliavin differentiable, and
for any $p'\geq 1$,
\begin{align}
\sup_{x\in\mR^d}\mE\left(\sup_{s\in[t,S]}\|DX_{t,s}(x)\|^{p'}_{\mH}\right)<+\infty,
\end{align}
where $D$ is the Malliavin derivative (cf. \cite{Nu}).
\item For any $f\in C^1_b(\mR^d)$, we have the following derivative formula: for Lebesgue-almost all $x\in\mR^d$,
\begin{align}
\nabla\mE f(X_{t,s}(x))=\frac{1}{s-t}\mE\left(f(X_{t,s}(x))\int^s_t\sigma^{-1}_r(X_{t,r}(x))\nabla X_{t,r}(x)\dif W_r\right),\label{EU05}
\end{align}
where $\sigma^{-1}$ is the inverse matrix of $\sigma$.
\item Assume that $b'\in\mL^q_p(T,S)$ with the same $p,q$ as 
in the assumptions. Let $X^{b,\sigma}_{t,s}(x)$ and $X^{b',\sigma}_{t,s}(x)$ be
the solutions to (\ref{SDE00}) associated with $b$ and $b'$ respectively. Then
\begin{align}
\sup_{x\in\mR^d}\mE\left(\sup_{s\in[t,S]}|X^{b,\sigma}_{t,s}(x)-X^{b',\sigma}_{t,s}(x)|^2\right)\leq C\|b-b'\|^2_{\mL^q_p(t,S)},\label{EU55}
\end{align}
where $C=C\big(d,p,q,K,\alpha,\|b\|_{\mL^q_p(t,S)},\|b'\|_{\mL^q_p(t,S)},\|\nabla\sigma\|_{\mL^q_p(t,S)}\big)$.
\end{enumerate}
\et
\br
Conclusions {\bf (A)} and {\bf (B)} are not really new and they are contained in \cite{Kr-Ro, Fe-Fl, Zh2}.
Conclusions {\bf (C)}, {\bf (D)} and {\bf (E)} seem to be new.
Our proofs are based on Zvonkin's transformation (cf. \cite{Zv}) and some results from the theory of PDEs.
The global $L^p$-integrability of the coefficients plays a crucial role in our argument.
It should be noticed that when $\sigma_t(x)=\sigma_t$ and $b_t(x)$ are bounded,
{\bf (A)}, {\bf (B)} and  {\bf (C)} were studied in \cite{Me} and \cite{Mo-Ni-Pr}
by using different arguments.
Moreover, unlike \cite{Zh0} and \cite{Zh2}, there is no explosion time problem here since
we are assuming global integrability conditions on $\sigma$ and $b$, see Lemma \ref{Le51} (4) below.
\er
\br
The stability estimate (\ref{EU55}) could be used to study
numerical solutions of SDEs with singular drifts. For example, let us consider the following SDE:
$$
\dif X_t=1_A(X_t)\dif t+\dif W_t,\ \ X_0=x,
$$
where $A$ is a bounded open subset of $\mR^d$. Let $b_n(x)=1_A*\varrho_n(x)$ be the mollifying approximation.
By (\ref{EU55}),
the solution $X^n_t$ of the above SDE corresponding to $b_n$ converges to $X_t$ in $L^2$.
Next, we can approximate $X^n_t$ by Euler's scheme.
In this way, one can give a numerical approximation for solutions of singular SDEs.
We plan to pursue this in a future project.
We would also like to mention that the
derivative formula \eqref{EU05} could be used in the computation of Greeks for pay-off functions in mathematical finance (cf. \cite{Ma-Th}).
\er
In the remainder of this section, we present an application of the above theorem to
incompressible Navier-Stokes equations. This application is actually one of the motivations of the present paper.
Consider the following classical Navier-Stokes equation in $\mR^3$:
$$
\p_t u=\nu\Delta u-(u\cdot\nabla) u+\nabla p,\ \ \div u=0,\ \ u_0=\varphi,
$$
where $u$ is the velocity field, $\nu$ is the viscosity constant and $p$ is the pressure of the fluid, $\varphi$ is the initial velocity with vanishing divergence.
In \cite{Co-Iy}, Constantin and Iyer provided a probabilistic representation to the above NSE as follows:
\begin{align}\label{NS22}
\left\{
\begin{aligned}
X_t(x)&=x+\int^t_0 u_s(X_s(x))\dif s+\sqrt{2\nu}W_t,\\
u_t(x)&=\bP\mE[\nabla^{\mathrm{t}} X^{-1}_t\cdot \varphi(X^{-1}_t)](x),
\end{aligned}
\right.
\end{align}
where $X^{-1}_t(x)$ denotes the inverse flow of $x\mapsto X_t(x)$, $\nabla^{\mathrm{t}} X^{-1}_t$ is the transpose of the Jacobian matrix,
and $\bP=\mI-\nabla(-\Delta)^{-1}\div$ is Leray's projection onto the space of all divergence free vector fields.
Let $\omega=\mathrm{curl}(u)=\nabla\times u$ be the vorticity.
Then the second equation in \eqref{NS22} can be written as
\begin{align}
\omega_t(x)=\mE[(\nabla X^{-1}_t(x))^{-1}\cdot \omega_0(X^{-1}_t(x))],\ \ \omega_0=\nabla\times \varphi,\label{NS33}
\end{align}
where $(\nabla X^{-1}_t(x))^{-1}$ stands for the inverse matrix of $\nabla X^{-1}_t(x)$.
In this case, the velocity $u$ can be recovered from $\omega$ by Biot-Savart's law (cf. \cite{Ma-Be}):
\begin{align}
u_t(x)=\int_{\mR^3}K_3(x-y)\omega_t(y)\dif y=:\bK\omega_t(x),\label{NS1}
\end{align}
where
$$
K_3(x) h=\frac{1}{4\pi}\frac{x\times h}{|x|^3},\  \ x,h\in\mR^3.
$$
In other words, we have the following stochastic representation to vorticity:
\begin{align}\label{NS2}
\left\{
\begin{aligned}
X_t(x)&=x+\int^t_0 \bK\omega_s(X_s(x))\dif s+\sqrt{2\nu}W_t,\\
\omega_t(x)&=\mE[(\nabla X^{-1}_t(x))^{-1}\cdot \omega_0(X^{-1}_t(x))].
\end{aligned}
\right.
\end{align}
Now if we substitute (\ref{NS33}) and (\ref{NS1}) into (\ref{NS2}), then we obtain the following equation:
$$
X_t(x)=x+\tilde\mE \int^t_0\!\!\!\int_{\mR^3}[K_3(X_s(x)-y)\nabla^{-1} \tilde X^{-1}_s(y)\cdot \omega_0(\tilde X^{-1}_s(y))]\dif y\dif s+\sqrt{2\nu}W_t,
$$
where the random field $\{\tilde X_t(y)\}_{y\in\mR^d}$ is an independent copy of $\{X_t(x)\}_{x\in\mR^d}$, and
$\tilde\mE$ denotes the expectation with respect to $(\tilde X_t)$ given $(X_t)$.
By the change of variables $\tilde X^{-1}_t(y)=x'$ and noticing that
$$
\det \nabla\tilde X_t(x')=1, \ \ (\nabla\tilde X^{-1}_t(\tilde X_t(x')))^{-1} =\nabla\tilde X_s(x'),
$$
we further have
$$
X_t(x)=x+\tilde\mE \int^t_0\!\!\!\int_{\mR^3}[K_3(X_s(x)-\tilde X_s(x'))\nabla\tilde X_s(x')\cdot \omega_0(x')]\dif x'\dif s+\sqrt{2\nu}W_t.
$$
This is simply the random vortex method for Navier-Stokes equations studied in \cite[Chapter 6]{Ma-Be}.

Recently, in \cite{Zh1} and \cite{Zh3}, we studied a backward analogue
of the stochastic representation (\ref{NS22}), that is, for $\nu>0$ and $t\leq s\leq 0$,
\begin{align}\label{NS3}
\left\{
\begin{aligned}
X_{t,s}(x)&=x+\int^s_t u_r(X_{t,r}(x))\dif r+\sqrt{2\nu}(W_s-W_t),\\
u_t(x)&=\bP\mE[\nabla^{\mathrm{t}}X_{t,0}\cdot\varphi(X_{t,0})](x).
\end{aligned}
\right.
\end{align}
The advantage of this representation is that the inverse of stochastic flow $x\mapsto X_{t,0}(x)$ does not appear.
In this case, $u_t(x)$ solves the following backward Navier-Stokes equation:
$$
\p_t u+\nu\Delta u-(u\cdot\nabla) u+\nabla p=0,\ \ \div u=0,\ \ u_0=\varphi,
$$
Using Theorem \ref{Main}, we have the following local well-posedness to the stochastic system (\ref{NS3}).
\bt\label{Th12}
For any $p>d$ and divergence free $\varphi\in \mW^1_p(\mR^d;\mR^d)$,
there exist a time $T=T(p,d,\nu,\|\varphi\|_{\mW^1_p})<0$ and a unique pair $(u,X)$ with $u\in L^\infty([T,0]; \mW^1_p)$
solving the stochastic system (\ref{NS3}).
\et

This paper is organized as follows: In Section 2, we recall some well-known results and give some preliminaries about the Sobolev differentiabilities of random vector fields.
In Section 3, we study a class of parabolic partial differential equations with time dependent coefficients and give some necessary estimates.
In Section 4, we prove some Krylov type  and Khasminskii type estimates.
In Section 5, we prove our main Theorem \ref{Main} for SDE (\ref{SDE00}) with $b=0$.
In Section 6, we prove Theorem \ref{Main}. In Section 7, we prove Theorem \ref{Th12} by using Theorem \ref{Main} and a fixed point argument.

Throughout this paper, we use the following convention: $C$ with or without subscripts will denote a positive constant, whose value may change in different places, and
whose dependence on the parameters can be traced from the calculations.

\section{Prelimiaries}
We first introduce some spaces and notations for later use. For $p,q\in[1,\infty]$ and $T<S$,
we denote by $\mL^q_p(T,S)$ the space of
all real-valued Borel functions on $[T,S]\times\mR^d$ with norm
$$
\|f\|_{\mL^q_p(T,S)}:=\left(\int^S_T\left(\int_{\mR^d}|f(t,x)|^p\dif x\right)^{\frac{q}{p}}\right)^{\frac{1}{q}}<+\infty.
$$
For $m\in\mN$ and $p\geq 1$, let $\mW^m_p=\mW^m_p(\mR^d)$ be the usual Sobolev space over $\mR^d$ with norm
$$
\|f\|_{\mW^m_p}:=\sum_{k=0}^m\|\nabla^k f\|_p<+\infty,
$$
where $\nabla^k$ denotes the $k$-order gradient operator, and $\|\cdot\|_p$ is the usual $L^p$-norm. For $\beta\geq 0$, let
$\mH^\beta_p:=(I-\Delta)^{-\frac{\beta}{2}}(L^p)$ be the usual Bessel potential space with norm (cf. \cite{St, Tr})
$$
\|f\|_{\mH^\beta_p}:=\|(I-\Delta)^{\frac{\beta}{2}}f\|_p.
$$
Notice that for $m\in\mN$ and $p>1$,
$$
\|f\|_{\mH^m_p}\asymp \|f\|_{\mW^m_p},
$$
where $\asymp$ means that the two sides are comparable up to a positive constant.
Moreover, let $\sC^\beta$ be the usual H\"older space with finite norm
$$
\|f\|_{\sC^\beta}:=\sum_{k=0}^{[\beta]}\|\nabla^k f\|_\infty+\sup_{x\not=y}\frac{|\nabla^{[\beta]}f(x)-\nabla f^{[\beta]} f(y)|}{|x-y|^{\beta-[\beta]}}<\infty,
$$
where $[\beta]$ is the integer part of $\beta$. By Sobolev's embedding theorem, we have
\begin{align}
\|f\|_{\sC^\delta}\leq C\|f\|_{\mH^\beta_p}, \ \ \beta-\delta>d/p, \ \delta\geq 0.\label{SOB}
\end{align}
In this paper we shall also use the following Banach space:
$$
\mW^{2,q}_p(T,S):=L^q(T,S;\mW^2_p)\cap \mW^{1,q}([T,S]; L^p).
$$

Let $f$ be a locally integrable function
on $\mR^d$. The Hardy-Littlewood maximal function is defined by
$$
\cM f(x):=\sup_{0<r<\infty}\frac{1}{|B_r|}\int_{B_r}f(x+y)\dif y,
$$
where $B_r:=\{x\in\mR^d: |x|<r\}$.
We recall the following result (cf. \cite[Appendix A]{Cr-De-Le}).
\bl\label{Le2}
(i) There exists a constant $C_d>0$ such that for all $f\in \mW^1_1(\mR^d)$ and Lebesgue-almost all $x,y\in \mR^d$,
\begin{align}
|f(x)-f(y)|\leq C_d |x-y|(\cM|\nabla f|(x)+\cM|\nabla f|(y)).\label{ES2}
\end{align}
(ii) For any $p>1$, there exists a constant $C_{d,p}>0$ such that for all $f\in L^p(\mR^d)$,
\begin{align}
\|\cM f\|_p\leq C_{d,p}\|f\|_p.\label{Es30}
\end{align}
\el
For $p>1$, let $\sV_p$ be the set of all continuous random fields $X:\mR^d\times\Omega\to\mR^d$ with
\begin{align}\label{LU1}
\|X\|_{\sV_p}:=\|X(0)\|_{L^p_\omega}+\|\nabla X\|_{L^\infty_x(L^p_\omega)}<\infty,
\end{align}
where $\nabla X$ denotes the generalized Jacobian matrix, and
$$
L^p_\omega:=L^p(\Omega),\ L^\infty_x(L^p_\omega):=L^\infty(\mR^d;L^p(\Omega)).
$$
Let $\sV^0_p\subset\sV_p$ be the set of random fields satisfying the additional condition
\begin{align}
\int_{\mR^d}\mE f(X(x))\dif x=\int_{\mR^d}f(x)\dif x.\label{CN1}
\end{align}
\br
The continuity assumption of $x\mapsto X(x)$ in the definition of $\sV_p$ is purely technical for $p>d$.
In fact, if $X\in\sV_p$ for $p>d$, then by Sobolev's embedding theorem, $x\mapsto X(x)$ always has a continuous version.
Condition (\ref{CN1}) means that $x\mapsto X(x)$ preserves the volume in the sense of mean values.
In the sequel, we also use the following notation:
$$
\sV_{\infty-}:=\cap_{p>1}\sV_p,\ \ \sV^0_{\infty-}:=\cap_{p>1}\sV^0_p,\ \ L^\infty_x(L^{\infty-}_\omega):=\cap_{p>1}L^\infty_x(L^p_\omega).
$$
\er
Let $\varrho:\mR^d\to[0,1]$ be a smooth function with support in $B_1$ and
$\int\varrho\dif x=1$. For $n\in\mN$, define a family of mollifiers $\varrho_n(x)$ as follows:
\begin{align}
\varrho_n(x):=n^{d}\varrho(nx),\ x\in\mR^d.\label{Rho}
\end{align}
For $X\in\sV_p$, define
\begin{align}
X_n(x):=\varrho_n*X(x)=\int_{\mR^d}X(x-y)\varrho_n(y)\dif y.\label{R1}
\end{align}
Clearly, by Jensen's inequality we have
\begin{align}
\sup_{x\in\mR^d}\mE|\nabla X_n(x)|^p\leq\mathrm{ess.}\sup_{x\in\mR^d}\mE|\nabla X(x)|^p=\|\nabla X\|^p_{L^\infty_x(L^p_\omega)}.\label{R2}
\end{align}
\bl
Let $p>1$. For any $X\in\sV_p$, we have
\begin{align}
\mE|X(x)-X(y)|^p\leq |x-y|^p\|\nabla X\|^p_{L^\infty_x(L^p_\omega)},\ \ \forall x,y\in\mR^d.\label{Le71}
\end{align}
\el
\begin{proof}
Let $X_n$ be defined by (\ref{R1}). By Fatou's lemma and (\ref{R2}), we have for all $x,y\in\mR^d$,
\begin{align*}
&\mE|X(x)-X(y)|^p\leq\varliminf_{n\to\infty}\mE|X_n(x)-X_n(y)|^p\\
&\leq|x-y|^p\varliminf_{n\to\infty}\int^1_0\mE|\nabla X_n(x+\theta(y-x))|^p\dif\theta\\
&\leq|x-y|^p\sup_{x\in\mR^d}\mE|\nabla X_n(x)|^p\leq|x-y|^p\|\nabla X\|^p_{L^\infty_x(L^p_\omega)},
\end{align*}
where we have used the continuity of $x\mapsto X(x)$ in the first inequality.
\end{proof}

\bl\label{Le24}
For any $p>1$, let $\{X_n,n\in\mN\}\subset\sV_p$ be a bounded sequence and $X(x)$ a continuous random field.
If, for each $x\in\mR^d$, $X_n(x)$ converges to $X(x)$ in probability,
then $X\in\sV_p$ and
$$
\|\nabla X\|_{L^\infty_x(L^p_\omega)}\leq\sup_n\|\nabla X_n\|_{L^\infty_x(L^p_\omega)}.
$$
Moreover, for some subsequence $n_k$, $\nabla X_{n_k}$ weakly converges to $\nabla X$ as random variables in $L^p(\Omega\times B_R;\mM^d)$ for any $R\in\mN$,
where $B_R=\{x: |x|<R\}$.
\el
\begin{proof}
Recall the definition of  $\sV_p$.
Since $\sup_n\|X_n(0)\|_{L^p_\omega}<\infty$, by (\ref{R2}) and (\ref{Le71}), we have for any $R>0$,
\begin{align}
\sup_n\int_{B_R}(\mE|X_n(x)|^p+\mE|\nabla X_n(x)|^p)\dif x<\infty.\label{R3}
\end{align}
This means that $\{X_n(\cdot), n\in\mN\}$ is bounded in $L^p(\Omega;\mW^1_p(B_R))$,
where $\mW^1_p(B_R)$ is the first-order Sobolev space over $B_R$.
Since $L^p(\Omega;\mW^1_p(B_R))$ is weakly compact, by a diagonal argument, there exist a subsequence $n_k$ and a random field
$\tilde X\in \cap_{R\in\mN}L^p(\Omega;\mW^1_p(B_R))$ such that for any $R\in\mN$,
\begin{align}
X_{n_k}(x)\to \tilde X(x)\mbox{ weakly in $L^p(\Omega;\mW^1_p(B_R))$}.\label{EU2}
\end{align}
In particular, for any $Z\in C^\infty_0(\mR^d; \mR^d)$ and $\xi\in L^\infty(\Omega)$, we have
$$
\lim_{k\to\infty}\mE\int_{\mR^d}\<X_{n_k}(x), Z(x)\xi \>_{\mR^d}\dif x=\mE\int_{\mR^d}\<\tilde X(x), Z(x)\xi \>_{\mR^d}\dif x.
$$
Since for each $x\in\mR^d$, $X_n(x)$ converges to $X(x)$ in probability, by (\ref{R3}) and the dominated convergence theorem, we also have
$$
\lim_{k\to\infty}\mE\int_{\mR^d}\<X_{n_k}(x), Z(x)\xi \>_{\mR^d}\dif x=\mE\int_{\mR^d}\<X(x), Z(x)\xi \>_{\mR^d}\dif x.
$$
Thus, for all $Z\in C^\infty_0(\mR^d; \mR^d)$ and $\xi\in L^\infty(\Omega)$,
$$
\mE\int_{\mR^d}\<X(x), Z(x)\xi \>_{\mR^d}\dif x=\mE\int_{\mR^d}\<\tilde X(x), Z(x)\xi \>_{\mR^d}\dif x,
$$
which implies that $X(x,\omega)=\tilde X(x,\omega)$ for $\dif x\times P(\dif \omega)$-almost all $(x,\omega)$. In particular, for almost all $\omega$,
$x\mapsto X(x,\omega)$ is Sobolev differentiable, and by (\ref{EU2}),
$\nabla X_{n_k}$ weakly converges to $\nabla X$ as random variables in $L^p(\Omega\times B_R;\mM^d)$ for each $R\in\mN$.

Now, let $\sV^\infty_c$ be the set of all $\mM^d$-valued smooth random fields with compact supports and bounded derivatives. Let $p_*=p/(p-1)$.
Since the dual space of $L^1(\mR^d; L^{p_*}(\Omega))$ is $L^\infty(\mR^d; L^{p}(\Omega))$
and $\sV^\infty_c$ is dense in $L^1(\mR^d; L^{p_*}(\Omega))$, we have
\begin{align*}
\|\nabla X\|_{L^\infty_x(L^p_\omega)}&=\sup_{U\in\sV_c^\infty;\|U\|_{L^1(L^{p_*})}\leq 1}\left|\int_{\mR^d}\mE\<\nabla X(x), U(x)\>_{\mM^d}\dif x\right|\\
&=\sup_{U\in\sV_c^\infty;\|U\|_{L^1(L^{p_*})}\leq 1}\left|\mE\left(\int_{\mR^d}\<X(x), \div U(x)\>_{\mR^d}\dif x\right)\right|\\
&=\sup_{U\in\sV_c^\infty;\|U\|_{L^1(L^{p_*})}\leq 1}\lim_{n\to\infty}\left|\mE\left(\int_{\mR^d}\<X_n(x), \div U(x)\>_{\mR^d}\dif x\right)\right|\\
&=\sup_{U\in\sV_c^\infty;\|U\|_{L^1(L^{p_*})}\leq 1}\lim_{n\to\infty}\left|\mE\left(\int_{\mR^d}\<\nabla X_n(x), U(x)\>_{\mM^d}\dif x\right)\right|\\
&\leq\sup_{n\in\mN}\sup_{U\in\sV_c^\infty;\|U\|_{L^1(L^{p_*})}\leq 1}\left|\mE\left(\int_{\mR^d}\<\nabla X_n(x), U(x)\>_{\mM^d}\dif x\right)\right|
=\sup_{n\in\mN}\|\nabla X_n\|_{L^\infty_x(L^p_\omega)}.
\end{align*}
The proof is complete.
\end{proof}

\bp\label{Pr74}
Let $p_1,p_2,p_3\in(1,\infty)$ with $\frac{1}{p_3}=\frac{1}{p_1}+\frac{1}{p_2}$.
If $X\in\sV_{p_1}$ and $Y\in\sV_{p_2}$ are two independent random fields, then we have $X\circ Y\in\sV_{p_3}$ and
\begin{align}
\|\nabla(X\circ Y)\|_{L^\infty_x(L^{p_3}_\omega)}\leq\|\nabla X\|_{L^\infty_x(L^{p_1}_\omega)}\|\nabla Y\|_{L^\infty_x(L^{p_2}_\omega)}.\label{R5}
\end{align}
Moreover, if for each $x\in\mR^d$, $\omega\mapsto X(x,\omega), Y(x,\omega)$ are Malliavin differentiable and
$$
\sup_{x\in\mR^d}\mE\|DX(x)\|_\mH^{p_1}<\infty,\ \sup_{x\in\mR^d}\mE\|DY(x)\|_\mH^{p_2}<\infty,
$$
then $X\circ Y(x)$ is also Malliavin differentiable and
\begin{align}
\sup_{x\in\mR^d}\mE\|D(X\circ Y(x))\|_\mH^{p_3}<\infty.\label{R6}
\end{align}
\ep
\begin{proof}
Let $X_n$ be defined by (\ref{R1}).  By (\ref{Le71}), we have
\begin{align*}
\sup_{x\in\mR^d}\mE|X_n(x)-X(x)|^{p_1}
&\leq\sup_{x\in\mR^d}\mE\int_{\mR^d}|X(x-y)-X(x)|^{p_1}\varrho_n(y)\dif y\\
&\leq\|\nabla X\|^{p_1}_{L^\infty_x(L^{p_1}_\omega)}\int_{\mR^d}|y|^{p_1}\rho_n(y)\dif y\leq\|\nabla X\|^{p_1}_{L^\infty_x(L^{p_1}_\omega)}/n^{p_1}.
\end{align*}
Since  $(X_n(x), X(x))_{x\in\mR^d}$ and $(Y_n(x), Y(x))_{x\in\mR^d}$ are independent,
we have for each $x\in\mR^d$,
\begin{align*}
\mE|X_n\circ Y(x)-X\circ Y(x)|^{p_1}&=\mE\left(\mE|X_n(y)-X(y)|^{p_1}|_{y=Y(x)}\right)\\
&\leq\sup_y\mE|X_n(y)-X(y)|^{p_1}\leq\|\nabla X\|^{p_1}_{L^\infty_x(L^{p_1}_\omega)}/n^{p_1}
\end{align*}
and
\begin{align*}
\|X_n\circ Y_n(x)-&X_n\circ Y(x)\|_{L^{p_3}_\omega}\leq\left\||Y_n(x)-Y(x)|\int^1_0|\nabla X_n|(Y_n(x)+\theta(Y(x)-Y_n(x)))\dif\theta\right\|_{L^{p_3}_\omega}\\
&\leq\|Y_n(x)-Y(x)\|_{L^{p_2}_\omega}\sup_x\|\nabla X_n(x)\|_{L^{p_1}_\omega}
\leq\|\nabla X\|_{L^\infty_x(L^{p_1}_\omega)}\|\nabla Y\|_{L^\infty_x(L^{p_2}_\omega)}/n.
\end{align*}
Since $p_3\leq p_1$, we thus have
\begin{align}
\lim_{n\to\infty}\sup_{x\in\mR^d}\mE|X_n\circ Y_n(x)-X\circ Y(x)|^{p_3}=0.\label{R7}
\end{align}
On the other hand, by the chain rule and H\"older's inequality, we have
\begin{align*}
\|\nabla (X_n\circ Y_n)\|_{L^\infty_x(L^{p_3}_\omega)}&\leq \sup_{x\in\mR^d}\left[(\mE|(\nabla X_n)\circ Y_n(x)|^{p_1})^{1/p_1}
(\mE|\nabla Y_n(x)|^{p_2})^{1/p_2}\right]\\
&\leq \|\nabla X_n\|_{L^\infty_x(L^{p_1}_\omega)}\|\nabla Y_n\|_{L^\infty_x(L^{p_2}_\omega)}
\leq \|\nabla X\|_{L^\infty_x(L^{p_1}_\omega)}\|\nabla Y\|_{L^\infty_x(L^{p_2}_\omega)},
\end{align*}
which, together with (\ref{R7}) and by Lemma \ref{Le24}, yields (\ref{R5}).

Similarly, by the chain rule,
$$
D(X_n\circ Y_n(x))=(DX_n)\circ Y_n(x)+\nabla X_n\circ Y_n(x)\cdot DY_n(x),
$$
and since $(DX_n(x), \nabla X_n(x))_{x\in\mR^d}$ and $(Y_n(x))_{x\in\mR^d}$ are independent, as above, we have
\begin{align*}
\|D(X_n\circ Y_n)\|_{L^\infty_x(L^{p_3}_\omega)}&\leq
\|(DX_n)\circ Y_n\|_{L^\infty_x(L^{p_3}_\omega)}+\|\nabla X_n\circ Y_n\cdot DY_n\|_{L^\infty_x(L^{p_3}_\omega)}\\
&\leq\|(DX_n)\circ Y_n\|_{L^\infty_x(L^{p_1}_\omega)}+\|\nabla X_n\circ Y_n\|_{L^\infty_x(L^{p_1}_\omega)}\|DY_n\|_{L^\infty_x(L^{p_2}_\omega)}\\
&\leq\|DX_n\|_{L^\infty_x(L^{p_1}_\omega)}+\|\nabla X_n\|_{L^\infty_x(L^{p_1}_\omega)}\|DY_n\|_{L^\infty_x(L^{p_2}_\omega)}\\
&\leq\|DX\|_{L^\infty_x(L^{p_1}_\omega)}+\|\nabla X\|_{L^\infty_x(L^{p_1}_\omega)}\|DY\|_{L^\infty_x(L^{p_2}_\omega)},
\end{align*}
which, together with (\ref{R7}) and by \cite[p.79, Lemma 1.5.3]{Nu}, yields (\ref{R6}).
\end{proof}

\section{A study of PDE $\p_t u+L^\sigma_tu+f=0$}

In the remainder of this paper, we shall fix $T<S$ with $S-T\leq 1$.
Suppose that $\sigma:[T,S]\times\mR^d\to\mM^d$ is a bounded Borel function.
Let us consider the following backward PDE:
\begin{align}
\p_tu+L^{\sigma}_tu+f=0,\ \ u(S)=0,\label{PDE}
\end{align}
where $f:[T,S]\times\mR^d\to\mR$ is a measurable function and
\begin{align}
L^{\sigma}_tu(x):=\tfrac{1}{2}\sigma^{ik}_t(x)\sigma^{jk}_t(x)\p_i\p_ju(x).\label{LL2}
\end{align}
Here and in the rest of this paper, we use the convention that repeated indices in a product
will be summed automatically.
The aim of this section is to prove
\bt\label{Th11}
Assume that $\sigma$ satisfies {\bf (H$^\alpha_K$)}. Let $p\in(1,\infty)$.
For any $f\in \mL^p_p(T,S)$, there exists a unique solution $u\in \mW^{2,p}_p(T,S)$ to (\ref{PDE}) with
\begin{align}
\|u\|_{\mL^p_p(T,S)}+\|\p_tu\|_{\mL^p_p(T,S)}+\|\nabla^2_xu\|_{\mL^p_p(T,S)}\leq C\|f\|_{\mL^p_p(T,S)},\label{Es1}
\end{align}
where $C=C(d,\alpha,K,p)>0$.
Furthermore, if $p,q\in(1,\infty)$ and $f\in \mL^p_p(T,S)\cap\mL^q_p(T,S)$,
then for any $\beta\in[0,2)$ and $\gamma>1$ with
$\frac{2}{q}+\frac{d}{p}<2-\beta+\frac{d}{\gamma}$,
\begin{align}\label{Es2}
\|u(t)\|_{\mH^\beta_{\gamma}}\leq C(S-t)^{\frac{2-\beta}{2}-\frac{d}{2p}-\frac{1}{q}+\frac{d}{2\gamma}}\|f\|_{\mL^q_p(t,S)},
\end{align}
where $C=C(d,\alpha,K,p,q,\gamma,\beta)$ is independent of $t\in[T,S]$.
\et

We first prove the a priori estimate \eqref{Es1}.
\bl
For any $p\in(1,\infty)$ and $f\in \mL^p_p(T,S)$, let $u\in \mW^{2,p}_p(T,S)$ satisfy (\ref{PDE}).
If $\sigma$ satisfies {\bf (H$^\alpha_K$)}, then \eqref{Es1} holds for some $C=C(d,\alpha,K,p)>0$.
In particular, the uniqueness holds for \eqref{PDE} in the class of $u\in\mW^{2,p}_p(T,S)$.
\el
\begin{proof}
We use the freezing coefficient argument (cf. \cite[Chapter 1]{Kr0}) and divide the proof into four steps.
\\
{\bf (1)} In this step, we first assume $\sigma_t(x)=\sigma_t$ does not depend on $x$.
For $f\in L^p(\mR^d)$, define
\begin{align}
\cT_{t,s}f(x):=\mE f\left(x+\int^s_t\sigma_r\dif W_r\right)=\int_{\mR^d}f(y)\rho(t,x;s,y)\dif y,\label{EG1}
\end{align}
where
$$
\rho(t,x;s,y)=\frac{\e^{-\<A_{t,s}^{-1}(x-y), x-y\>/2}}{\sqrt{(2\pi)^d\det(A_{t,s})}},\ \ A_{t,s}:=\int^s_t\sigma^\mathrm{t}_r\sigma_r\dif r.
$$
In this case, the unique solution of \eqref{PDE} is explicitly given by
\begin{align}
u(t,x)=\int^S_t\cT_{t,s}f(s,x)\dif s.\label{Es22}
\end{align}
By \cite[Theorem 1.1]{Kr1}, for any $p,q\in(1,\infty)$, there exists a constant $C_0=C_0(d,K,p,q)>0$ such that
\begin{align}
\left(\int^S_T\left\|\nabla^2_x\int^S_t\cT_{t,s}f(s,\cdot)\dif s\right\|^q_{p}\dif t\right)^{1/q}\leq C_0\|f\|_{\mL^q_p(T,S)}.\label{Es11}
\end{align}
{\bf (2)} Next, we assume that for some $x_0\in\mR^d$,
\begin{align}
\|\sigma_t(x)-\sigma_t(x_0)\|\leq \tfrac{1}{2C_0 K},\label{LK1}
\end{align}
where $C_0$ is the constant in \eqref{Es11} and $K$ is the constant in {\bf (H$^\alpha_K$)}.
In this case, we may write
$$
\p_t u+L^{\sigma_\cdot(x_0)}_t u+g=0,\mbox{ where } \ g:=L^{\sigma}_t u-L^{\sigma_\cdot(x_0)}_t u+f.
$$
Note that by the definition of $L^{\sigma}_t$ and \eqref{LK1},
$$
\|g\|_{\mL^q_p(T,S)}\leq \tfrac{1}{2C_0}\|\nabla^2_x u\|_{\mL^q_p(T,S)}+\|f\|_{\mL^q_p(T,S)}.
$$
Thus, by \eqref{Es22} and \eqref{Es11}, we have
\begin{align*}
\|\nabla^2_x u\|_{\mL^q_p(T,S)}\leq C_0\|g\|_{\mL^q_p(T,S)}
\leq \tfrac{1}{2}\|\nabla^2_x u\|_{\mL^q_p(T,S)}+C_0\|f\|_{\mL^q_p(T,S)},
\end{align*}
which in turn gives
$$
\|\nabla^2_x u\|_{\mL^q_p(T,S)}\leq 2C_0\|f\|_{\mL^q_p(T,S)}.
$$
{\bf (3)} Let $\zeta:\mR^d\to[0,1]$ be a smooth function with $\zeta(x)=1$ for $|x|\leq 1$ and $\zeta(x)=0$ for $|x|\geq 2$.
Fix a small constant $\delta$ whose value will be determined below.
For fixed $z\in\mR^d$, set
$$
\zeta^\delta_z(x):=\zeta((x-z)/\delta).
$$
It is easy to see that for $j=0,1,2,$
\begin{align}
\int_{\mR^d}|\nabla^j_x\zeta^\delta_z(x)|^p\dif z=\delta^{d-jp}\int_{\mR^d}|\nabla^j\zeta(z)|^p\dif z>0.\label{EK5}
\end{align}
Multiplying both sides of (\ref{PDE}) by $\zeta^\delta_z$, we obtain
\begin{align}\label{EK2}
\p_t(u\zeta^\delta_z)+L^{\sigma}_t (u\zeta^\delta_z)+g^\delta_z=0,
\end{align}
where
$$
g^\delta_z:=L^{\sigma}_t(u\zeta^\delta_z)-(L^\sigma_t u)\zeta^\delta_z+f\zeta^\delta_z.
$$
Define
$$
\tilde\sigma_t(x):=\sigma_t\big((x-z)\zeta^{2\delta}_z(x)+z\big).
$$
Since $\zeta^\delta_z(x)=1$ for $|x-z|\leq \delta$ and $\zeta^\delta_z(x)=0$ for $|x-z|>2\delta$, we have
\begin{align}\label{EK3}
L^{\sigma}_t (u\zeta^\delta_z)=L^{\tilde\sigma}_t (u\zeta^\delta_z).
\end{align}
Notice that by {\bf (H$^\alpha_K$)},
$$
\|\tilde\sigma_t(x)-\tilde\sigma_t(z)\|\leq K|(x-z)\zeta^{2\delta}_z|^\alpha\leq K|4\delta|^\alpha,
$$
and
$$
\|g^\delta_z\|_{\mL^q_p}\leq K^2\||\nabla_x u|\cdot|\nabla_x\zeta^\delta_z|\|_{\mL^q_p}+K^2\||u|\cdot|\nabla^2_x\zeta^\delta_z|\|_{\mL^q_p}+\|f\zeta^\delta_z\|_{\mL^q_p}.
$$
Letting $\delta$ be small enough, by \eqref{EK2}, \eqref{EK3} and step (2), we have
\begin{align}
\|\nabla^2_x(u\zeta^\delta_z)\|_{\mL^q_p(t,S)}&\leq 2C_0\|g^\delta_z\|_{\mL^q_p(t,S)}\leq 2C_0K^2\||\nabla_x u|\cdot|\nabla_x\zeta^\delta_z|\|_{\mL^q_p(t,S)}\no\\
&+2C_0K^2\||u|\cdot|\nabla^2_x\zeta^\delta_z|\|_{\mL^q_p(t,S)}+2C_0\|f\zeta^\delta_z\|_{\mL^q_p(t,S)}.\label{EK4}
\end{align}
{\bf (4)} If $p=q$, then integrating both sides of \eqref{EK4} with respect to $z$,
and using \eqref{EK5} and Fubini's theorem, we obtain
$$
\int_{\mR^d}\|\nabla^2_x(u\zeta^\delta_z)\|^p_{\mL^p_p(t,S)}\dif z\leq C\Big(\|\nabla_x u\|^p_{\mL^p_p(t,S)}+\|u\|^p_{\mL^p_p(t,S)}+\|f\|^p_{\mL^p_p(t,S)}\Big).
$$
Hence, by \eqref{EK5} again, $\|\nabla u\|_p\leq C\|\nabla^2 u\|^{\frac{1}{2}}_p\|u\|^{\frac{1}{2}}_p$ and Young's inequality, we have
\begin{align*}
&\|\nabla^2_xu\|^p_{\mL^p_p(t,S)}=\int_{\mR^d}\|\nabla^2_xu\cdot\zeta^\delta_z\|^p_{\mL^p_p(t,S)}\dif z\\
&\quad\leq C\Big(\|\nabla_x u\|^p_{\mL^p_p(t,S)}+\|u\|^p_{\mL^p_p(t,S)}+\|f\|^p_{\mL^p_p(t,S)}\Big)\\
&\quad\leq \tfrac{1}{2}\|\nabla^2_x u\|^p_{\mL^p_p(t,S)}+C\Big(\|u\|^p_{\mL^p_p(t,S)}+\|f\|^p_{\mL^p_p(t,S)}\Big).
\end{align*}
Thus, for some $C=C(d,\alpha,K,p)>0$,
\begin{align}\label{LU2}
\|\nabla^2_xu\|^p_{\mL^p_p(t,S)}\leq C\Big(\|u\|^p_{\mL^p_p(t,S)}+\|f\|^p_{\mL^p_p(t,S)}\Big),
\end{align}
which together with \eqref{PDE} gives
$$
\|u(t)\|^p_{p}\leq  C\|u\|^p_{\mL^p_p(t,S)}+C\|f\|^p_{\mL^p_p(T,S)}=C\int^S_t\|u(s)\|^p_{p}\dif s+C\|f\|^p_{\mL^p_p(T,S)}.
$$
By Gronwall's inequality, \eqref{LU2} and \eqref{PDE}, we obtain \eqref{Es1}.
\end{proof}
\br
In the above proof, the reason we required $p=q$ was due to the use of Fubini's theorem.
In the case $p\not=q$,
it seems that we can not use the freezing coefficient argument to obtain the
a priori estimate \eqref{Es1}
since in general it is not true that for some $\gamma\in[1,\infty]$,
$$\int_{\mR^d}\|f\cdot\zeta^\delta_z\|^\gamma_{\mL^q_p(t,S)}\dif z\asymp \|f\|^\gamma_{\mL^q_p(t,S)}.$$
We leave \eqref{Es1} for $p\not=q$ as an open problem.
\er
Next we show the existence of a solution to \eqref{PDE}  in $\mW^{2,p}_p(T,S)$
and \eqref{Es2} by using mollifying and weak convergence arguments. For this purpose we assume $\sigma$
satisfies {\bf (H$^\alpha_K$)} and for some $\alpha'\in(0,1)$ and $K'>0$,
\begin{align}\label{Tim}
\|\sigma_t(x)-\sigma_s(x)\|\leq K'|t-s|^{\alpha'}.
\end{align}
Under {\bf (H$^\alpha_K$)} and \eqref{Tim}, it is a classical fact that
the operator $\p_t+L^\sigma_t$ has a fundamental solution $\rho(t,x;s,y)$
(see e.g. \cite[Chapter IV]{La-So-Ur} or \cite[Chapter 1]{Fr}), i.e., for any $f\in C_b(\mR^d)$,
the function
$$
\cT_{t,s}f(x):=\int_{\mR^d}f(y)\rho(t,x;s,y)\dif y
$$
satisfies that for all $(t,x)\in[T,S]\times\mR^d$,
\begin{align}\label{Heat}
\p_t\cT_{t,s}f(x)+L^\sigma_t\cT_{t,s}f(x)=0,\ \ \lim_{t\uparrow s}\cT_{t,s}f(x)=f(x).
\end{align}
Furthermore, for all $x,y\in\mR^d$ and $T\leq t<s\leq S$ (see \cite[p.376, (13.1)]{La-So-Ur}),
\begin{align}
|\nabla^j_x\rho(t,x;s,y)|\leq C_j (s-t)^{-\frac{j}{2}}(2(s-t)^{-\frac{d}{2}}\e^{-\kappa_j|x-y|^2/(2(s-t))},\ \ j=0,1,2,\label{Es4}
\end{align}
where $C_j, \kappa_j>0$ only depend on $\alpha, K$ and $d$.

Here is an easy corollary of (\ref{Es4}).
\bl\label{Le23}
For any $p,\gamma\in(1,\infty)$ and $\beta\in[0,2)$, there exists a constant $C=C(d,\alpha,K,p,\gamma,\beta)>0$
such that for all $f\in L^p(\mR^d)$ and $T\leq t<s\leq S$,
\begin{align}
\|\cT_{t,s}f\|_{\mH^\beta_\gamma}\leq C(s-t)^{-\frac{\beta}{2}-\frac{d}{2p}+\frac{d}{2\gamma}}\|f\|_p.\label{Es6}
\end{align}
\el
\begin{proof}
By the heat kernel estimate (\ref{Es4}), we have for all $p\in[1,\infty]$,
$$
\|\nabla^j\cT_{t,s}f\|_p\leq C(s-t)^{-\frac{j}{2}}\|f\|_p,\ \ j=0,1,2.
$$
By Gagliardo-Nirenberg's and complex interpolation inequalities (cf. \cite[Theorem 2.1]{Ta}), we have
$$
\|\cT_{t,s}f\|_{\mH^\beta_\gamma}\leq C\|\nabla^2\cT_{t,s}f\|^{\frac{\beta}{2}+\frac{d}{2p}-\frac{d}{2\gamma}}_{p}
\|\cT_{t,s}f\|^{\frac{2-\beta}{2}-\frac{d}{2p}+\frac{d}{2\gamma}}_{p}
\leq C(s-t)^{-\frac{\beta}{2}-\frac{d}{2p}+\frac{d}{2\gamma}}\|f\|_p,
$$
which gives (\ref{Es6}).
\end{proof}

Let $f\in C([T,S]; \mW^{2}_p)$ and define
$$
u(t,x):=\int^S_t\cT_{t,s}f(s,x)\dif s.
$$
By \eqref{Heat}, it is easy to see that $u\in\mW^{2,p}_p(T, S)$ satisfies \eqref{PDE}. Moreover, for any $p,q,\gamma\in(1,\infty)$ and
$\beta\in[0,2)$  with $\frac{2}{q}+\frac{d}{p}<2-\beta+\frac{d}{\gamma}$, by \eqref{Es6} and H\"older's inequality, we have
\begin{align}\label{LU4}
\begin{split}
\|u(t)\|_{\mH^\beta_\gamma}&\leq\int^S_t\|\cT_{t,s}f(s)\|_{\mH^\beta_\gamma}\dif s\leq C\int^S_t(s-t)^{-\frac{\beta}{2}-\frac{d}{2p}+\frac{d}{2\gamma}}\|f(s)\|_p\dif s\\
&\leq C\left(\int^S_t(s-t)^{-\frac{\beta q^*}{2}-\frac{dq^*}{2p}+\frac{dq^*}{2\gamma}}\dif s\right)^{\frac{1}{q^*}}\|f\|_{\mL^q_p(t,S)}\\
&\leq C(S-t)^{\frac{2-\beta}{2}-\frac{d}{2p}-\frac{1}{q}+\frac{d}{2\gamma}}\|f\|_{\mL^q_p(t,S)},
\end{split}
\end{align}
where $q^*:=\frac{q}{q-1}$ and $C=C(d,\alpha,K,p,q,\gamma,\beta)>0$.

\

Now we are ready to give

\begin{proof}[Proof of Theorem \ref{Th11}]
Let $\varrho$ be a nonnegative smooth function in $\mR^{d+1}$ with support in $\{x\in\mR^{d+1}: |x|\leq 1\}$ and
$\int_{\mR^{d+1}}\varrho(t,x)\dif t\dif x=1$. Set $\varrho_n(t,x):=n^{d+1}\varrho(nt,nx)$ and extend $u(s)$ to $\mR$ by setting
$u(s,\cdot)=0$ for $s\notin[T,S]$. Define
\begin{align}
\sigma_n:=\sigma*\varrho_n,\ \ f_n:=f*\varrho_n.\label{BP6}
\end{align}
Let $u_n$ solve the following equation
\begin{align}\label{APP}
\p_tu_n+L^{\sigma_n}_tu_n+f_n=0,\ \ u_n(S)=0.
\end{align}
By \eqref{Es1} and \eqref{LU4}, we have the following uniform estimate:
\begin{align}
\|u_n\|_{\mL^p_p(T,S)}+\|\p_tu_n\|_{\mL^p_p(T,S)}+\|\nabla^2_xu_n\|_{\mL^p_p(T,S)}\leq C\|f\|_{\mL^p_p(T,S)},\label{Es81}
\end{align}
and for any $\beta\in[0,2)$ and $\gamma, q>1$ with $\frac{2}{q}+\frac{d}{p}<2-\beta+\frac{d}{\gamma}$,
\begin{align}\label{Es82}
\|u_n(t)\|_{\mH^\beta_\gamma}\leq C(S-t)^{\frac{2-\beta}{2}-\frac{d}{2p}-\frac{1}{q}+\frac{d}{2\gamma}}\|f\|_{\mL^q_p(t,S)},
\end{align}
where the constant $C$ only depends on $d,\alpha,K, p,q,\gamma,\beta$.

By \eqref{Es81} and the weak compactness of $\mW^{2,p}_p(T,S)$,
there exist a subsequence still denoted by $u_n$ and a function $u\in\mW^{2,p}_p(T,S)$ with $u(S)=0$ such that
$u_n$ weakly converges to $u$. By taking weak limits of \eqref{APP}, one sees that $u$ satisfies \eqref{PDE}. Indeed,
for any $\varphi\in C^\infty_0((T,S)\times\mR^{d})$, we have
\begin{align*}
\left|\int^S_T\!\!\!\int_{\mR^{d}}(L^{\sigma_m}_tu_n-L^{\sigma}_tu_n)\varphi\dif t\dif x\right|
&\leq C\left(\int^S_T\|\sigma_m(t)-\sigma(t)\|_\infty\|\nabla^2_xu_n\|_p\dif t\right)\\
&\leq C\left(\int^S_T\|\sigma_m(t)-\sigma(t)\|^{\frac{p}{p-1}}_\infty\dif t\right)^{\frac{p-1}{p}}\|\nabla^2_xu_n\|_{\mL^p_p(T,S)},
\end{align*}
which,  by \eqref{Es81}, converges to zero as $m\to\infty$  uniformly in $n$. On the other hand, for fixed $m$, since $u_n$ weakly converges to $u$, we have
$$
\int^S_T\!\!\!\int_{\mR^{d}}(L^{\sigma_m}_tu_n-L^{\sigma_m}_tu)\varphi\dif t\dif x\to 0,\ \mbox{ as $n\to\infty$.}
$$
Hence,
$$
\int^S_T\!\!\!\int_{\mR^{d}}(L^{\sigma_n}_tu_n-L^{\sigma}_tu)\varphi\dif t\dif x\to 0,\ \mbox{ as $n\to\infty$.}
$$
Similarly, for any $\varphi\in C^\infty_0((T,S)\times\mR^{d})$, we have
\begin{align*}
\int^S_T\!\!\!\int_{\mR^{d}}(\p_tu_n)\varphi\dif t\dif x&=-\int^S_T\!\!\!\int_{\mR^{d}}u_n\p_t\varphi\dif t\dif x\\
&\to-\int^S_T\!\!\!\int_{\mR^{d}}u\p_t\varphi\dif t\dif x=\int^S_T\!\!\!\int_{\mR^{d}}\p_t u\varphi\dif t\dif x
\end{align*}
as $n\to\infty$, and by the property of convolutions,
$$
\lim_{n\to\infty}\|f_n-f\|_{\mL^p_p(T,S)}=0.
$$
Moreover, as in the proof of Lemma \ref{Le24}, by \eqref{Es82}
we get \eqref{Es2}.
\end{proof}

\section{Krylov type and Khasminskii type estimates}
The following Krylov estimate was proved in \cite[Theorem 2.1]{Zh2}.
Since we need more explicit dependence on $s-t$, for the reader's convenience, we reproduce the proof here.
\bt\label{Krylov}
Assume that $\sigma$ satisfies {\bf (H$^\alpha_K$)} and $q,p\in(1,\infty)$ with $\frac{d}{p}+\frac{2}{q}<2$. Let $0<S-T\leq 1$.
For any $s\in[T,S]$ and $x\in\mR^d$, let $X_{T, s}(x)$ solve SDE (\ref{SDE00}) with $b=0$. For any $\delta\in(0,1-\frac{d}{2p}-\frac{1}{q})$,
there exists a positive constant $C=C(K,\alpha,d,p,q,\delta)$
such that for all $f\in\mL^q_p(T,S)$, $T\leq t\leq s\leq S$  and $x\in\mR^d$,
\begin{align}
\mE\left(\int^s_t f(r,X_{T,r}(x))\dif r\Bigg|_{\sF_t}\right)\leq C(s-t)^{\delta}\|f\|_{\mL^q_p(T,S)},\label{HG1}
\end{align}
where $\sF_t:=\sigma\{W_s: s\leq t\}$.
\et
\begin{proof}
Let $p'=2d$. Since $\mL^{p'}_{p'}(T, S)\cap\mL^q_p(T, S)$ is dense in $\mL^q_p(T,S)$, it suffices to prove (\ref{HG1}) for
$$
f\in\mL^{p'}_{p'}(T,S)\cap\mL^q_p(T,S).
$$
Fix $s\in[T,S]$. By Theorem \ref{Th11}, there exists a unique solution $u\in\mW^{2,p'}_{p'}(T,s)$ to the following backward PDE:
$$
\p_t u+L^\sigma_tu+f=0,\ t\in[T,s], \ u(s,x)=0,
$$
so that for all $t\in[T,s]$,
$$
\|u\|_{\mL^{p'}_{p'}(t,s)}+\|\nabla^2u\|_{\mL^{p'}_{p'}(t,s)}\leq C\|f\|_{\mL^{p'}_{p'}(t,s)}.
$$
Moreover, by (\ref{Es2}) and \eqref{SOB}, for any $\delta\in(0,1-\frac{d}{2p}-\frac{1}{q})$, we have
\begin{align}
\sup_{r\in[t,s]}\|u(r)\|_\infty\leq C(s-t)^{\delta}\|f\|_{\mL^q_p(t,s)}, \ \ \forall t\in[T,s].\label{BP7}
\end{align}
Let $\varrho_n$ be the same mollifiers as in the proof of Theorem \ref{Th11}. Define
\begin{align}
u_n(t,x):=u*\varrho_n(t,x),\ \ f_n(t,x):=-[\p_t u_n(t,x)+L^\sigma_tu_n(t,x)].\label{HG2}
\end{align}
Then we have
\begin{align*}
\|f_n-f\|_{\mL^{p'}_{p'}(t,s)}&\leq\|\p_t(u_n-u)\|_{\mL^{p'}_{p'}(t,s)}+K\|\nabla^2(u_n-u)\|_{\mL^{p'}_{p'}(t,s)}\\
&\leq \|\p_tu*\varrho_n-\p_tu\|_{\mL^{p'}_{p'}(t,s)}+K\|\nabla^2 u*\varrho_n-\nabla^2u\|_{\mL^{p'}_{p'}(t,s)}\\
&\leq\|f*\varrho_n-f\|_{\mL^{p'}_{p'}(t,s)}+2K\|\nabla^2 u*\varrho_n-\nabla^2u\|_{\mL^{p'}_{p'}(t,s)},
\end{align*}
which converges to zero as $n\to\infty$ by the property of convolutions.
So, by the classical Krylov estimate (cf. \cite[Lemma 5.1]{Kr} or \cite[Lemma 3.1]{Gy-Ma}), we have
\begin{align}
\lim_{n\to\infty}\mE\left(\int^{s}_t |f_n(r,X_{T,r})-f(r,X_{T,r})|\dif r\right)
\leq C\lim_{n\to\infty}\|f_n-f\|_{\mL^{p'}_{p'}(t,s)}=0.\label{BP5}
\end{align}
Now applying It\^o's formula to $u_n(t,x)$ and using (\ref{HG2}), we get that
for any $T\leq t\leq s\leq S$,
$$
u_n(s,X_{T,s})=u_n(t,X_{T,t})-\int^s_t f_n(r,X_{T,r})\dif r+\int^s_t\p_i u_n(r,X_{T,r})\sigma^{ik}_r(X_{T,r})\dif W^k_r.
$$
Since
$$
\sup_{s,x}|\p_i u_n(s,x)|\leq C_n,
$$
by Doob's optional theorem, we have
$$
\mE\left[\int^s_t\p_i u_n(r,X_{T,r})\sigma^{ik}_r(X_{T,r})\dif W^k_r\Bigg|_{\sF_t}\right]=0.
$$
Hence,
\begin{align*}
\mE\left(\int^s_t f_n(r,X_{T,r})\dif r\Bigg|_{\sF_t}\right)
&=\mE\Bigg[(u_n(t,X_{T,t})-u_n(s,X_{T,s}))\Big|_{\sF_t}\Bigg]\no\\
&\leq 2\sup_{(r,x)\in[t,s]\times\mR^d}|u_n(r,x)|\leq 2\sup_{r\in[t,s]}\|u(r)\|_\infty\\
&\leq C(s-t)^{\delta}\|f\|_{\mL^q_p(T,S)},
\end{align*}
where the last step is due to (\ref{BP7}).
Combining this with (\ref{BP5}) we arrive at the desired conclusion.
\end{proof}
We also need the following Khasminskii type estimate (cf. \cite[Lemma 1.1]{Po}).
\bl\label{Le35}
Let $(\xi(t))_{t\in [S, T]}$, $(\eta(t))_{t\in [S, T]}$ and $(\beta(t))_{t\in [S, T]}$ be
three real-valued measurable $\sF_t$-adapted processes, 
and $(\eta(t))_{t\in [S, T]}$ and $(\alpha(t))_{t\in [S, T]}$ be two $\mR^d$-valued measurable $\sF_t$-adapted processes.
Suppose that there exist $c_0>0$ and $\delta\in(0,1)$ such that for any $T\leq t\leq s\leq S$
\begin{align}
\mE\left(\int^s_t[|\beta(r)|+|\alpha(r)|^2]\dif r\Big|\sF_t\right)\leq c_0(s-t)^\delta,\label{C1}
\end{align}
and that
$$
\xi(s)=\xi(T)+\int^s_T\zeta(r)\dif r+\int^s_T\eta(r)\dif W_r+\int^s_T\xi(r)\beta(r)\dif r+\int^s_T\xi(r)\alpha(r)\dif W_r.
$$
Then for any $p>0$ and $\gamma_1,\gamma_2,\gamma_3>1$, we have
\begin{align}
\mE\left(\sup_{s\in[T,S]}\xi^+(s)^p\right)&\leq C\left(\|\xi^+(T)^p\|_{\gamma_1}
+\left\|\left(\int^S_T\zeta^+(r)\dif r\right)^p\right\|_{\gamma_2}
+\left\|\left(\int^S_T|\eta(r)|^2\dif r\right)^{\frac{p}{2}}\right\|_{\gamma_3}\right),\label{HG8}
\end{align}
where $a^+=\max\{0,a\}$, $C=C(c_0,\delta,p,\gamma_i)>0$ and $\|\cdot\|_\gamma$ denotes the norm in $L^\gamma(\Omega)$.
\el
\begin{proof}
Write
$$
M(s):=\exp\left\{\int^s_T\alpha(r)\dif W_r-\frac{1}{2}\int^s_T|\alpha(r)|^2\dif r+\int^s_T\beta(r)\dif r\right\}.
$$
By It\^o's formula, one sees that
\begin{align}
\xi(s)=M(s)\left\{\xi(T)+\int^s_TM^{-1}(r)\Big(\eta(r)\dif W_r+[\zeta(r)-\<\alpha(r),\eta(r)\>]\dif r\Big)\right\}.\label{C2}
\end{align}
By (\ref{C1}) and the Khasminskii estimate (cf. \cite[Lemma 1.1]{Po}),
we have for any $p\geq 1$,
$$
\mE\exp\left\{p\int^S_T|\alpha(r)|^2\dif r+p\int^S_T|\beta(r)|\dif r\right\}\leq C=C(c_0,\beta,p)<\infty,
$$
which implies that for any $p\in\mR$,
$$
s\mapsto\exp\left\{p\int^s_T\alpha(r)\dif W_r-\frac{p^2}{2}\int^s_T|\alpha(r)|^2\dif r\right\}
$$
is an exponential martingale. Thus, by H\"older's inequality and Doob's maximal inequality, 
we have that for any $p\in\mR$,
$$
\mE\left(\sup_{s\in[T,S]}|M(s)|^p\right)\leq C=C(c_0,\delta,p)<\infty.
$$
The desired estimate follows by (\ref{C2}), H\"older  and Burkholder's inequalities.
\end{proof}

\section{SDEs without drifts}

In this section, we consider the following SDE:
\begin{align}
\dif X_{t,s}=\sigma_s(X_{t,s})\dif W_s,\ X_{t,t}=x, \ s\geq t,\label{SDE0}
\end{align}
where $\sigma:[T,S]\times\mR^d\to\mM^d$ satisfies {\bf (H$^\alpha_K$)}.
It is well-known that, under {\bf (H$^\alpha_K$)}, (\ref{SDE0})
is well-posed in the sense of Stroock-Varadhan's martingale solutions (cf. \cite[p187, Theorem 7.2.1]{St-Va}).
Indeed, H\"older's continuity can be replaced with the weaker condition that $\sigma$ is 
uniformly continuous in $x$ with respect to $t$.
Moreover, $\{X_{t,s}(x)\}$ defines a family of time non-homogeneous Markov processes.
The aim of this section is to prove Theorem \ref{Main} for SDE (\ref{SDE0}). More precisely, we want to prove
\bt\label{Th0}
Assume that $\sigma$ satisfies {\bf (H$^\alpha_K$)} and  that
for some $q, p\in(1,\infty)$ with $\frac{d}{p}+\frac{2}{q}<1$,
$$
\nabla\sigma_t\in \mL^q_p(T,S).
$$
Then we have the following conclusions:
\begin{enumerate}[{\bf (a)}]
\item For any $(t,x)\in[T,S]\times\mR^d$, there is a unique strong solution denoted 
by $X_{t,s}(x)$ or $X^\sigma_{t,s}(x)$ to (\ref{SDE0}), which has a jointly continuous version with respect to $s,x$.
\item For each $s\geq t$ and almost all $\omega$, 
$x\mapsto X_{t,s}(x,\omega)$ is weakly differentiable.  Let $\nabla X_{t,s}(x)$
be the Jacobian matrix and $J_{t,s}(x)$ solve the following linear matrix-valued SDE:
\begin{align}
J_{t,s}(x)=\mI+\int^s_t\nabla\sigma_r(X_{t,r}(x))J_{t,r}(x)\dif W_r.\label{EU6}
\end{align}
Then $J_{t,s}(x)=\nabla X_{t,s}(x)$ a.s. for Lebesgue almost all $x\in\mR^d$, and for any $p'\geq 1$,
\begin{align}
\sup_{x\in\mR^d}\mE\left(\sup_{s\in[t,S]}|J_{t,s}(x)|^{p'}\right)\leq C=C(p,q,d,K,\alpha,p',\|\nabla\sigma\|_{\mL^q_p(T,S)}),\label{EU4}
\end{align}
where the constant $C$ is increasing with respect to $\|\nabla\sigma\|_{\mL^q_p(T,S)}$.
\item For each $s\geq t$ and $x\in\mR^d$, the random variable $\omega\mapsto X_{t,s}(x,\omega)$ is Malliavin differentiable, and
for any $p'\geq 1$,
\begin{align}
\sup_{x\in\mR^d}\mE\left(\sup_{s\in[t,S]}\|DX_{t,s}(x)\|^{p'}_{\mH}\right)\leq C=C(p,q,d,K,\alpha,p',\|\nabla\sigma\|_{\mL^q_p(T,S)}).\label{EU94}
\end{align}
Moreover, for any adapted vector field $h$ with $\mE\int^S_T|\dot h(r)|^2\dif r<\infty$, the Malliavin derivative $D_hX_{t,s}(x)$ along $h$ satisfies the following linear SDE:
\begin{align}
D_hX_{t,s}(x)=\int^s_t\nabla\sigma_r(X_{t,r}(x))D_hX_{t,r}(x)\dif W_r+\int^s_t\sigma_r(X_{t,r}(x))\dot h(r)\dif r.\label{EU96}
\end{align}
\item For any $f\in C^1_b(\mR^d)$, we have the following formula: for Lebesgue almost all $x\in\mR^d$,
\begin{align}
\nabla\mE f(X_{t,s}(x))=\frac{1}{s-t}\mE\left(f(X_{t,s}(x))\int^s_t\sigma^{-1}_r(X_{t,r}(x))\nabla X_{t,r}(x)\dif W_r\right).\label{EU5}
\end{align}
\item Assume that $\sigma'$ 
also satisfies the assumptions of the theorem with the same $K,\alpha$ and $p,q$.
Let $X^\sigma_{t,s}(x)$ and $X^{\sigma'}_{t,s}(x)$ be the solutions to (\ref{SDE0})
associated with $\sigma$ and $\sigma'$ respectively. Then 
$$
\sup_{x\in\mR^d}\mE\left(\sup_{s\in[t,S]}|X^\sigma_{t,s}(x)-X^{\sigma'}_{t,s}(x)|^2\right)\leq C(S-t)^\delta\|\sigma-\sigma'\|^2_{\mL^q_p(t,S)},
$$
provided $\|\sigma-\sigma'\|^2_{\mL^q_p(t,S)}<\infty$, where $\delta\in(0,1)$ only depends on $p,q,d$.
\end{enumerate}
\et
\subsection{Some a priori estimates} In this subsection, we assume that $\sigma$ satisfies {\bf (H$^\alpha_K$)} and
$$
\sup_{t,x}|\nabla^j\sigma_t(x)|<\infty,\ \forall j\in\mN.
$$
In this case, it is well-known that the unique solution $X^\sigma_{t,s}(x)$ (or simply denoted by
$X_{t,s}$) of (\ref{SDE0}) 
forms a $C^\infty$-diffeomorphism flow  (cf. \cite[p.312, Theorem 39]{Pr}).
Let $J_{t,s}:=\nabla X_{t,s}$ be the Jacobian matrix, and $DX_{t,s}$ the Malliavin derivative of $X_{t,s}$ with respect to sample paths.
Then we have (cf. \cite[p.312, Theorem 39]{Pr})
\begin{align}
J_{t,s}=\mI+\int^s_t\nabla\sigma_r(X_{t,r})J_{t,r}\dif W_r,\label{HG7}
\end{align}
and for any $h\in\mH$,
\begin{align}
D_hX_{t,s}=\int^s_t\nabla\sigma_r(X_{t,r})D_hX_{t,r}\dif W_r
+\int^s_t\sigma_r(X_{t,r})\dot h_r\dif r.\label{HG77}
\end{align}

We have the following a priori estimates.
\bp
Under the assumptions of Theorem \ref{Th0}, 
for any $p'\geq 1$, we have
\begin{align}
\sup_{x\in\mR^d}\mE\left(\sup_{s\in[t,S]}|\nabla X_{t,s}(x)|^{p'}\right)
+\sup_{x\in\mR^d}\mE\left(\sup_{s\in[t,S]}\|D X_{t,s}(x)\|^{p'}_{\mH}\right)\leq C,\label{EU01}
\end{align}
where the constant $C=C(K,\alpha,p,q,d,p',\|\nabla\sigma\|_{\mL^{q}_{p}(T,S)})$ is increasing with respect to $\|\nabla\sigma\|_{\mL^q_p(T,S)}$.
\ep
\begin{proof}
Without loss of generality, we assume $t=T$ and write $X_s:=X_{T,s}$ and $J_s:=J_{T,s}$.
\\
{\bf (1)} Let
\begin{align*}
\beta(r):=\|\nabla\sigma_r(X_r)J_r\|^2/|J_r|^2,\ \ \alpha(r):=2\<J_r, \nabla\sigma_r(X_r)J_r\>/|J_r|^2.
\end{align*}
Here we use the convention $\frac{0}{0}:=0$, i.e., if $|J_r|=0$, then $\beta(r)=\alpha(r)=0$. By (\ref{HG7}) and It\^o's formula, we have
$$
|J_s|^2=|J_T|^2+\int^s_T|J_r|^2\beta(r)\dif r+\int^s_T|J_r|^2\alpha(r)\dif W_r.
$$
Let $\delta\in(0,1-\frac{d}{p}-\frac{2}{q})$. By (\ref{HG1}), we have for any $T\leq t\leq s\leq S$,
\begin{align*}
&\mE\left(\int^s_t\Big[|\alpha(r)|^2+|\beta(r)|\Big]\dif r\Bigg|_{\sF_t}\right)
\leq5\mE\left(\int^s_t|\nabla\sigma_r(X_r)|^2\dif r\Bigg|_{\sF_t}\right)\\
&\qquad\leq  C(s-t)^\delta\||\nabla\sigma|^2\|_{\mL^{q/2}_{p/2}(T,S)}=C(s-t)^\delta\|\nabla\sigma\|^2_{\mL^{q}_{p}(T,S)},
\end{align*}
which in turn gives the first estimate in (\ref{EU01}) by (\ref{HG8}).
\\
\\
{\bf (2)} For $T\leq r\leq s\leq S$, let $J_{r,s}$ solve the following linear SDE:
$$
J_{r,s}=\mI+\int^s_r\nabla\sigma_{r'}(X_{r'})J_{r,r'}\dif W_{r'}.
$$
By (\ref{HG77}) and the variation of constants formula, we have
\begin{align}
D_hX_s=\int^s_T J_{r,s}\sigma_r(X_r)\dot h_r\dif r.\label{PL1}
\end{align}
Let $\Sigma_s^{ij}:=\<DX^i_s, DX^j_s\>_\mH$ be the Malliavin covariance matrix.
Then by (\ref{PL1}), we have
\begin{align}
\Sigma_s=\int^s_T J_{r,s}\sigma_r(X_r)(J_{r,s}\sigma_r(X_r))^{\mathrm{t}}\dif r.\label{TP1}
\end{align}
As in step {\bf (1)}, one can show that for any $p'\geq 1$,
\begin{align}
\sup_{r\in[T,S]}\mE\left(\sup_{s\in[r,S]}|J_{r,s}|^{p'}\right)\leq C.\label{TP2}
\end{align}
Thus, by (\ref{TP1}) and (\ref{TP2}) we have
\begin{align*}
\mE\left(\sup_{s\in[T,S]}|\Sigma_s|^{p'}\right)&\leq C\mE\left(\sup_{s\in[T,S]}\int^s_T |J_{r,s}|^{2p'}\dif r\right)
\leq C\mE\left(\int^S_T\sup_{s\in[r,S]}|J_{r,s}|^{2p'}\dif r\right)\leq C.
\end{align*}
The proof is now complete.
\end{proof}

\bl
Assume that $\sigma, \sigma':[T,S]\times\mR^d\to \mM^d$ satisfy {\bf (H$^\alpha_K$)}
with the same $K,\alpha$. If for some $p,q\in(2,\infty)$
with $\frac{d}{p}+\frac{2}{q}<1$,
$$
\nabla\sigma_t,\ \nabla\sigma'_t\in \mL^q_p(T,S),
$$
then there exists a constant $C=C(K,\alpha,p,d,q,\|\nabla\sigma\|_{\mL^q_p(T,S)}, 
\|\nabla\sigma'\|_{\mL^q_p(T,S)})>0$ such that
\begin{align}
\sup_{x\in\mR^d}\mE\left(\sup_{s\in[t,S]}|X^\sigma_{t,s}(x)-X^{\sigma'}_{t,s}(x)|^2\right)\leq C(S-t)^\delta\|\sigma-\sigma'\|^2_{\mL^q_p(t,S)},\label{HG5}
\end{align}
where $\delta\in(0,1)$ only depends on $p,q,d$. Moreover, for any $\gamma>1$ and $x\in\mR^d$,
\begin{align}
\mE\left(\sup_{s\in[t,S]}|\nabla X^\sigma_{t,s}(x)-\nabla X^{\sigma'}_{t,s}(x)|^2\right)\leq
C\left\|\int^S_t|\nabla\sigma_r(X^{\sigma}_{t,r}(x))-\nabla\sigma'_r(X^{\sigma'}_{t,r}(x))|^2\dif r\right\|_{L^\gamma(\Omega)}.\label{HG6}
\end{align}
\el
\begin{proof}
Without loss of generality, we assume $t=T$ and write $X^\sigma_s:=X^\sigma_{T,s}$.
\\
{\bf (1)} Set $Z_s:=X^\sigma_s-X^{\sigma'}_s$, then
$$
Z_s=\int^s_T\Big[\sigma_r(X^{\sigma}_r)-\sigma'_r(X^{\sigma'}_r)\Big]\dif W_r.
$$
By It\^o's formula, we have
\begin{align*}
|Z_s|^2&=\int^s_T\|\sigma(r,X^{\sigma}_r)-\sigma'_r(X^{\sigma'}_r)\|^2\dif r+2\int^s_T\Big[\sigma(r,X^{\sigma}_r)-\sigma'_r(X^{\sigma'}_r)\Big]^{\mathrm{t}}Z_r\dif W_r\\
&=\int^s_T\zeta(r)\dif r+\int^s_T\eta(r)\dif W_r+\int^s_T|Z_r|^2\beta(r)\dif r+\int^s_T|Z_r|^2\alpha(r)\dif W_r,
\end{align*}
where
\begin{align*}
\zeta(r)&:=\|\sigma_r(X^{\sigma}_r)-\sigma'_r(X^{\sigma'}_r)\|^2-2\|\sigma_r(X^{\sigma}_r)-\sigma_r(X^{\sigma'}_r)\|^2,\\
\eta(r)&:=2[\sigma(r,X^{\sigma'}_r)-\sigma'_r(X^{\sigma'}_r)]^{\mathrm{t}}Z_r,\\
\beta(r)&:=2\|\sigma_r(X^\sigma_r)-\sigma_r(X^{\sigma'}_r)\|^2/|Z_r|^2,\\
\alpha(r)&:=2[\sigma_r(X^\sigma_r)-\sigma_r(X^{\sigma'}_r)]^{\mathrm{t}}Z_r/|Z_r|^2.
\end{align*}
Here we have used the convention $\frac{0}{0}:=0$, i.e., if $|Z_r|=0$, then $\beta(r)=\alpha(r)=0$.

By Lemma \ref{Le2}, (\ref{HG1}) and (\ref{Es30}), we have that for any $T\leq t<s\leq S$,
\begin{align*}
\mE\left(\int^s_t\Big[|\beta(r)|+|\alpha(r)|^2\Big]\dif r\Big|\sF_t\right)&\leq
C\mE\left(\int^s_t\Big[\cM|\nabla\sigma_r|^2(X^\sigma_r)+\cM|\nabla\sigma_r|^2(X^{\sigma'}_r)\Big]\dif r\Big|\sF_t\right)\\
&\leq C(s-t)^{\delta}\|\cM|\nabla\sigma|^2\|_{\mL^{q/2}_{p/2}(T,S)}\\
&\leq C(s-t)^{\delta}\||\nabla\sigma|^2\|_{\mL^{q/2}_{p/2}(T,S)}\\
&=C(s-t)^{\delta}\|\nabla\sigma\|^2_{\mL^{q}_{p}(T,S)},
\end{align*}
where $\delta\in(0,1-\frac{d}{p}-\frac{2}{q})$, and that for any $\gamma\in(1,1/(2/q+d/p))$,
\begin{align}
\mE\left(\int^S_T\|\sigma_r(X^{\sigma'}_r)-\sigma'_r(X^{\sigma'}_r)\|^{2\gamma}\dif r\right)
&\leq C(S-T)^{\delta}\|\|\sigma-\sigma'\|^{2\gamma}\|_{\mL^{q/(2\gamma)}_{p/(2\gamma)}(T,S)}\no\\
&=C(S-T)^{\delta}\|\sigma-\sigma'\|^{2\gamma}_{\mL^{q}_{p}(T,S)},\label{HG4}
\end{align}
where $\delta\in(0,1-\frac{d\gamma}{p}-\frac{2\gamma}{q})$.
Since
$\zeta^+(r)\leq 2\|\sigma_r(X^{\sigma'}_r)-\sigma'_r(X^{\sigma'}_r)\|^2,$
using (\ref{HG8}) with $p=1$, $\gamma_2=\gamma$ and $\gamma_3=\frac{2\gamma}{\gamma+1}$ and by H\"older's inequality,
we obtain
\begin{align}
\mE\left(\sup_{s\in[T,S]}|Z_s|^{2}\right)&\leq C\left\|\left(\int^S_T|Z_r|^2\|\sigma_r(X^{\sigma'}_r)
-\sigma'_r(X^{\sigma'}_r)\|^2\dif r\right)^{\frac{1}{2}}\right\|_{L^{\gamma_3}(\Omega)}\no\\
&\quad+C\left\|\int^S_T\|\sigma_r(X^{\sigma'}_r)-\sigma'_r(X^{\sigma'}_r)\|^2\dif r\right\|_{L^{\gamma_2}(\Omega)}\no\\
&\leq C\left\|\sup_{r\in[T,S]}|Z_r|\right\|_{L^2(\Omega)}
\left\|\int^S_T\|\sigma_r(X^{\sigma'}_r)-\sigma'_r(X^{\sigma'}_r)\|^2\dif r\right\|^{\frac{1}{2}}_{L^\gamma(\Omega)}\no\\
&\quad+C\left\|\int^S_T\|\sigma_r(X^{\sigma'}_r)-\sigma'_r(X^{\sigma'}_r)\|^2\dif r\right\|_{L^\gamma(\Omega)}\no\\
&\leq \frac{1}{2}\left\|\sup_{r\in[T,S]}|Z_r|\right\|^2_{L^2(\Omega)}+C\left\|\int^S_T\|\sigma_r(X^{\sigma'}_r)
-\sigma'_r(X^{\sigma'}_r)\|^2\dif r\right\|_{L^\gamma(\Omega)},\label{HG9}
\end{align}
which, together with (\ref{HG4}), yields (\ref{HG5}).
\\
\\
{\bf (2)} Set $U_s:=J^\sigma_s-J^{\sigma'}_s$. Then by (\ref{HG7}), we have
$$
U_s=\int^s_T[\nabla\sigma_r(X^{\sigma}_r)J^{\sigma}_r-\nabla\sigma'_r(X^{\sigma'}_r)J^{\sigma'}_r]\dif W_r.
$$
By It\^o's formula, we have
\begin{align*}
|U_s|^2&=2\int^s_T\<U_r, [\nabla\sigma_r(X^{\sigma}_r)J^{\sigma}_r-\nabla\sigma'_r(X^{\sigma'}_r)J^{\sigma'}_r]\dif W_r\>\\
&\quad+\int^s_T\|\nabla\sigma_r(X^{\sigma}_r)J^{\sigma}_r-\nabla\sigma'_r(X^{\sigma'}_r)J^{\sigma'}_r\|^2\dif r.
\end{align*}
As in the proof of (\ref{HG9}),
and using (\ref{EU01}) and by H\"older's inequality, we obtain that for $\gamma'>\gamma>1$,
\begin{align*}
\mE\left(\sup_{s\in[0,S]}|U_s|^2\right)&\leq C\left\|\int^S_T
\|[\nabla\sigma_r(X^{\sigma}_r)-\nabla\sigma'_r(X^{\sigma'}_r)]J^{\sigma'}_r\|^2\dif r\right\|_{L^\gamma(\Omega)}\\
&\leq C\left\|\sup_{r\in[T,S]}|J^{\sigma'}_r|^2\int^S_T|\nabla\sigma_r(X^{\sigma}_r)-\nabla\sigma'_r(X^{\sigma'}_r)|^2\dif r\right\|_{L^\gamma(\Omega)}\\
&\leq C\left\|\int^S_T|\nabla\sigma_r(X^{\sigma}_r)-\nabla\sigma'_r(X^{\sigma'}_r)|^2\dif r\right\|_{L^{\gamma'}(\Omega)},
\end{align*}
which gives (\ref{HG6}) by changing $\gamma'$ to $\gamma$.
\end{proof}
\subsection{Proof of Theorem \ref{Th0}}
{\bf (a)} Under the assumptions, the pathwise uniqueness follows from {\bf (e)}. Since $\sigma$ is bounded and uniformly continuous in $x$ with respect to $t$,
the existence of a weak solution is classical (cf. \cite{St}). The existence of a strong solution then follows by Yamada-Watanabe's theorem (cf. \cite[p163, Theorem 1.1]{Ik-Wa}).
\\
\\
{\bf (b)} Define $\sigma^n_t(x):=\sigma_t*\varrho_n(x)$, where $\varrho_n$ is a mollifier in $\mR^d$.
Consider the following SDE:
$$
X^n_{t,s}(x)=x+\int^s_t\sigma^n_r(X^n_{t,r}(x))\dif W_r,\ \ s\geq t.
$$
Since $\sigma^n$ is uniformly bounded, it is easy to see that for any $p'>1$,
$$
\sup_n\mE\left(\sup_{s\in[t,S]}|X^n_{t,s}(x)|^{p'}\right)\leq C(1+|x|^{p'}).
$$
Moreover, by (\ref{EU01}) we have
$$
\sup_{n}\sup_{x\in\mR^d}\mE\left(\sup_{s\in[t,S]}|\nabla X^n_{t,s}(x)|^{p'}\right)<\infty,
$$
and by (\ref{HG5}),
\begin{align}
\lim_{n\to\infty}\sup_{x\in\mR^d}\mE\left(\sup_{s\in[t,S]}|X^n_{t,s}(x)-X_{t,s}(x)|^2\right)\leq C\lim_{n\to\infty}\|\sigma^n-\sigma\|^2_{\mL^q_p(t,S)}=0.\label{C3}
\end{align}
Thus, by Lemma \ref{Le24}, the random field $x\mapsto X_{t,s}(x,\omega)$ is weakly differentiable almost surely, and for some subsequence $n_k$ and any $R\in\mN$,
\begin{align}
\mbox{$\nabla X^{n_k}_{t,s}$ weakly converges to $\nabla X_{t,s}$ as random variables in $L^{p'}(\Omega\times B_R;\mM^d)$}.\label{UY4}
\end{align}
Let $J_{t,s}(x)$ be the solution of SDE (\ref{EU6}). We need to show that $\nabla X_{t,s}(x)=J_{t,s}(x)$.  As in the proof of (\ref{EU01}), we have
$$
\sup_{x\in\mR^d}\mE\left(\sup_{s\in[t,S]}|J_{t,s}(x)|^{p'}\right)\leq C.
$$
Moreover, letting $J^n_{t,s}(x):=\nabla X^n_{t,s}(x)$, by (\ref{HG6}) we have
\begin{align}\label{D1}
\begin{split}
&\mE\left(\sup_{s\in[t,S]}|J^n_{t,s}(x)-J_{t,s}(x)|^2\right)\\
&\qquad\leq C\left\|\int^S_t|\nabla\sigma^n_r(X^n_{t,r}(x))-\nabla\sigma_r(X_{t,r}(x))|^2\dif r\right\|_{L^\gamma(\Omega)}.
\end{split}
\end{align}
As in the proof of (\ref{HG4}), we have for $\gamma\in(1,1/(2/q+d/p))$,
\begin{align}
\sup_{x\in\mR^d}\left\|\int^S_t|\nabla\sigma^m_r(X^n_{t,r}(x))-\nabla\sigma_r(X^n_{t,r}(x))|^2\dif r\right\|_{L^\gamma(\Omega)}\leq
C\|\nabla\sigma^m-\nabla\sigma\|^2_{\mL^q_p(t,S)},\label{D2}
\end{align}
where $C$ is independent of $n$. On the other hand, for fixed $m$, by (\ref{C3}) we have
\begin{align}
\lim_{n\to\infty}\sup_{x\in\mR^d}\left\|\int^S_t|\nabla\sigma^m_r(X^n_{t,r}(x))-\nabla\sigma^m_r(X_{t,r}(x))|^2\dif r\right\|_{L^\gamma(\Omega)}=0.\label{D3}
\end{align}
Combining (\ref{D1})-(\ref{D3}), we obtain
\begin{align}
\lim_{n\to\infty}\sup_{x\in\mR^d}\mE\left(\sup_{s\in[t,S]}|J^n_{t,s}(x)-J_{t,s}(x)|^2\right)=0,\label{PL2}
\end{align}
which, together with (\ref{UY4}), implies $\nabla X_{t,s}=J_{t,s}$ a.e.
\\
\\
{\bf (c)} By (\ref{EU01}) again, we have for any $p'\geq 1$,
$$
\sup_n\sup_{x\in\mR^d}\mE\left(\sup_{s\in[t,S]}\|DX^n_{t,s}(x)\|^{p'}_\mH\right)\leq C,
$$
which, together with (\ref{C3}) and by \cite[p.79, Lemma 1.5.3]{Nu}, yields that $X_{t,s}(x)$ is Malliavin differentiable and (\ref{EU94}) holds.
Let $h$ be an adapted vector field with $\mE\int^S_T|\dot h(r)|^2\dif r<\infty$. Then we have
$$
D_hX^n_{t,s}=\int^s_t\nabla\sigma_r(X^n_{t,r})D_hX^n_{t,r}\dif W_r+\int^s_t\sigma^n_r(X^n_{t,r})\dot h_r\dif r.
$$
Let $Z^h_{t,s}$ solve
$$
Z^h_{t,s}=\int^s_t\nabla\sigma_r(X_{t,r})Z^h_{t,r}\dif W_r+\int^s_t\sigma_r(X_{t,r})\dot h_r\dif r.
$$
As above, one can show that $D_hX^n_{t,s}\to Z^h_{t,s}$ in $L^2(\Omega)$. Moreover, for some subsequence $n_k$, $D_h X^{n_k}_{t,s}$ also weakly
converges to $D_hX_{t,s}$ in $L^2(\Omega)$. Thus, $Z^h_{t,s}=D_hX_{t,s}$ satisfies equation (\ref{EU96}).
\\
\\
{\bf (d)} By the classical Bismut-Elworthy-Li's formula (cf. \cite{El-Li}), we have for any $f\in C^1_b(\mR^d)$,
$$
\nabla\mE f(X^n_{t,s}(x))=\frac{1}{s-t}\mE\left[f(X^n_{t,s}(x))\int^s_t[\sigma^n_r(X^n_{t,r}(x))]^{-1}\nabla X^n_{t,r}(x)\dif W_r\right].
$$
Using (\ref{C3}) and (\ref{PL2}), by taking limits on both sides of the 
above formula, we obtain (\ref{EU5}). A more direct way of proving (\ref{EU5})
is to use {\bf (b)} and {\bf (c)}. We give it as follows: For fixed $\v\in\mR^d$ and $T\leq t<s\leq S$, define an adapted Cameron-Martin vector field $h_\v$ by
$$
h_\v(s'):=\frac{1}{s-t}\int^{s'}_t[\sigma_r(X_{t,r})]^{-1}\nabla_\v X_{t,r}\dif r,\ \ s'\in[t,s],
$$
where $\nabla_\v X_{t,r}:=\<\nabla X_{t,r},\v\>_{\mR^d}=J_{t,r}\v$.
By (\ref{EU4}), we have
$$
\mE\int^s_t|\dot h_\v(r)|^2\dif r=\frac{1}{(s-t)^2}\mE\int^{s}_t|[\sigma_r(X_{t,r})]^{-1}\nabla_\v X_{t,r}|^2\dif r<\infty.
$$
Notice that by \eqref{EU96}, $D_{h_\v}X_{t,{s'}}$ satisfies
$$
D_{h_\v}X_{t,s'}=\int^{s'}_t\nabla\sigma_r(X_{t,r})D_{h_\v}X_{t,r}\dif W_r+\frac{1}{s-t}\int^{s'}_t\nabla_\v X_{t,r}\dif r,\ \ s'\in[t,s].
$$
By (\ref{EU6}) and the variation  of constants formula, we have
$$
D_{h_\v}X_{t,s}=\nabla_\v X_{t,s}=J_{t,r}\v.
$$
Hence, by the chain rule and the integration by parts formula in the Malliavin calculus (cf. \cite{Nu}), we obtain
\begin{align*}
\nabla_\v\mE f(X_{t,s})&=\mE[\nabla f(X_{t,s})\nabla_\v X_{t,s}]
=\mE[\nabla f(X_{t,s})D_{h_\v}X_{t,s}]=\mE[D_{h_\v}(f(X_{t,s}))]\\
&=\frac{1}{s-t}\mE\left(f(X_{t,s})\int^s_t[\sigma_r(X_{t,r})]^{-1}\nabla_\v X_{t,r}\dif W_r\right).
\end{align*}
\\
{\bf (e)} 
Using \eqref{C3} and taking limits in 
$$
\mE\left(\sup_{s\in[t,S]}|X^{\sigma_n}_{t,s}(x)-X^{\sigma'_n}_{t,s}(x)|^2\right)\leq C(S-t)^\delta\|\sigma_n-\sigma'_n\|^2_{\mL^q_p(t,S)},
$$
we immediately get the desired conclusion.
The proof is now complete.

\section{Proof of Theorem \ref{Main}}

In this section we assume that $\sigma$ satisfies {\bf (H$^\alpha_K$)}
and that  one of the following two conditions holds:
\begin{enumerate}[(i)]
\item $\sigma_t(x)=\sigma_t$ is independent of $x$ and for some $p,q\in(1,\infty)$ with $\tfrac{d}{p}+\tfrac{2}{q}<1$, $b\in \mL^q_p(T,S).$
\item $\nabla\sigma, b\in \mL^q_p(T,S)$ for some $q=p>d+2$.
\end{enumerate}

We first prove the following result.
\bt\label{Th101}
Under the above assumptions (i) or (ii), for any $f\in \mL^q_p(T,S)$,  
there exists a unique solution $u=u^b_f\in \mW^{2,q}_p(T,S)$ to 
\begin{align}
\p_tu+L^{\sigma}_tu+b\cdot\nabla u+f=0,\ \ u(S)=0,\label{PDE1}
\end{align}
satisfying
\begin{align}
\|u\|_{\mL^q_p(T,S)}+\|\nabla^2u\|_{\mL^q_p(T,S)}\leq C_1\exp\left\{C_1\|b\|^q_{\mL^q_p(T,S)}\right\}\|f\|_{\mL^q_p(T,S)},\label{KH1}
\end{align}
and for all $t\in[T,S]$,
\begin{align}
\|\nabla u(t)\|_{\sC^{\delta/2}}\leq C_1(S-T)^{\delta/3}\exp\left\{C_1(S-T)^{q\delta/3}\|b\|^q_{\mL^q_p(T,S)}\right\}\|f\|_{\mL^q_p(T,S)},\label{ER1}
\end{align}
where $\delta:=\frac{1}{2}-\frac{d}{2p}-\frac{1}{q}$ and $C_1=C_1(K,\alpha,p,q,d,\delta)>0$. 
Suppose that $b'$ also satisfies the assumptions of this theorem and $f'\in \mL^q_p(T,S)$. Let $u^b_f$ and $u^{b'}_{f'}$
be the solutions of \eqref{PDE1} associated with $b, f$ and $b', f'$ respectively. Then
\begin{align}
\begin{split}
&\sum_{j=0,1}\|\nabla^j u^b_f(t)-\nabla^j u^{b'}_{f'}(t)\|_\infty+\sum_{j=0,2}\|\nabla^j u^b_f-\nabla^j u^{b'}_{f'}\|_{\mL^q_p(T,S)}\\
&\qquad\leq C_2\Big(\|f-f'\|_{\mL^q_p(T,S)}+\|b-b'\|_{\mL^q_p(T,S)}\Big),
\end{split}\label{ER11}
\end{align}
where $C_2=C_2(K,\alpha,p,q,d,\|b\|_{\mL^q_p(T,S)},\|b'\|_{\mL^q_p(T,S)}, \|f'\|_{\mL^q_p(T,S)})$.
\et
\begin{proof}
By standard Picard's iteration or a fixed point argument, we only need to prove the 
a priori estimates (\ref{KH1}), (\ref{ER1}) and (\ref{ER11}).
Letting $\delta:=\frac{1}{2}-\frac{d}{2p}-\frac{1}{q}$, by (\ref{Es2}), \eqref{SOB} with suitable choices of $\beta$ and $\gamma$, we have
\begin{align*}
\|\nabla u(t)\|^q_{\sC^{\delta/2}}&\leq C(S-T)^{q\delta/3}\int^S_t\|(b\cdot\nabla u)(s)+f(s)\|^q_p\dif s\\
&\leq C(S-T)^{q\delta/3}\int^S_t\Big[\|b(s)\|^q_p\|\nabla u(s)\|_\infty^q+\|f(s)\|^q_p\Big]\dif s,
\end{align*}
which, together with Gronwall's inequality, yields (\ref{ER1}).

On the other hand, in the case of (i), by (\ref{Es11}) and (\ref{ER1}), we have
\begin{align*}
&\|u\|_{\mL^q_p(T,S)}+\|\nabla^2u\|_{\mL^q_p(T,S)}\leq C\|(b\cdot\nabla u)+f\|_{\mL^q_p(T,S)}\\
&\quad\leq C\|b\|_{\mL^q_p(T,S)}\|\nabla u\|_\infty+C\|f\|_{\mL^q_p(T,S)}\\
&\quad\leq C\left(\|b\|_{\mL^q_p(T,S)}\exp\left\{C\|b\|^q_{\mL^q_p(T,S)}\right\}+1\right)\|f\|_{\mL^q_p(T,S)},
\end{align*}
which in turn gives (\ref{KH1}). In the case of (ii), by (\ref{Es1}) we still have (\ref{KH1}).

Moreover, if we let $w:=u^b_f-u^{b'}_{f'}$, then
\begin{align*}
\p_tw+L^\sigma_t w+b\cdot\nabla w+(b-b')\cdot\nabla u^{b'}_{f'}+f-f'=0,\ \ w(S)=0.
\end{align*}
As above, using  (\ref{Es2}), (\ref{SOB}) and \eqref{ER1}, and by Gronwall's inequality, we have
\begin{align*}
\|\nabla w\|_\infty&\leq C_1\exp\left\{C\Big(\|b\|^q_{\mL^q_p(T,S)}+\|b'\|^q_{\mL^q_p(T,S)}\Big)\right\}(\|f'\|_{\mL^q_p(T,S)}+1)\\
&\quad\times\Big(\|f-f'\|_{\mL^q_p(T,S)}+\|b-b'\|_{\mL^q_p(T,S)}\Big).
\end{align*}
The desired estimate (\ref{ER11}) follows by (\ref{Es2}),  \eqref{SOB} and (\ref{Es1}).
\end{proof}
Let $[t_0, s_0]\subset[T,S]$ be any subinterval. For $\ell=1,\cdots,d$, 
by Theorem \ref{Th101}, the following PDE
$$
\p_tu^\ell+L^\sigma_tu^\ell+b\cdot\nabla u^\ell+b^\ell=0,\ \ u^\ell_{s_0}(x)=0
$$
has a unique solution $u^\ell$.
Let
$$
\u_t(x):=\u^b_t(x):=(u^1_t(x),\cdots,u^d_t(x))
$$
and
\begin{align}
\Phi_t(x):=\Phi^b_t(x):=x+\u^b_t(x).\label{Def1}
\end{align}

We now prove the following Zvonkin transformation.
\bl\label{Le51}
Under (i) or (ii), for any $U>0$, there is a positive constant $\eps=\eps(K,\alpha,d,p,q,U)$ 
such that if $s_0-t_0\leq\eps$ and $\|b\|_{\mL^q_p(t_0,s_0)}\leq U$, then for each $t\in[t_0,s_0]$, $x\mapsto \Phi_t(x)$ is a $C^1$-diffeomorphism with
\begin{align}
\tfrac{1}{2}|x-y|\leq|\Phi_t(x)-\Phi_t(y)|\leq \tfrac{3}{2}|x-y|.\label{EQ1}
\end{align}
Moreover, letting $\delta:=\frac{1}{2}-\frac{d}{2p}-\frac{1}{q}>0$, we have the following conclusions:
\begin{enumerate}[{\bf (1)}]
\item $\|\nabla\Phi_t\|_\infty+\|\nabla\Phi^{-1}_t\|_\infty\leq \kappa$, where $\kappa$ is a universal constant.
\item $\|\nabla^2\Phi\|_{\mL^q_p(t_0,s_0)}+\|\nabla\Phi\|_{\sC^{\delta/2}}\leq C$, where $C$ only depends on $K,\alpha,p,q,d,\delta,U$.
\item Let $b'\in\mL^q_p(t_0,s_0)$ be another function with $\|b'\|_{\mL^q_p(t_0,s_0)}\leq U$. 
Let $\Phi^b$ and $\Phi^{b'}$ be associated with $b$ and $b'$ respectively. Then
we have
$$
\|\Phi^b-\Phi^{b'}\|_{\mL^\infty_\infty(t_0,s_0)}+\|\nabla\Phi^b-\nabla\Phi^{b'}\|_{\mL^q_p(t_0,s_0)}\leq C\|b-b'\|_{\mL^q_p(t_0,s_0)}.
$$
\item $X_{t_0,s}$ solves SDE (\ref{SDE00}) on $[t_0,s_0]$
if and only if $Y_{t_0,s}:=\Phi_s(X_{t_0,s})$ solves the following SDE:
\begin{align}
\dif Y_{t_0,s}=\Theta_s(Y_{t_0,s})\dif W_s,\ s\in[t_0,s_0],\ \ Y_{t_0,t_0}=\Phi_{t_0}(x),\label{SDE02}
\end{align}
where $\Theta_s(y):=[\nabla\Phi_s\cdot \sigma_s]\circ(\Phi^{-1}_s(y))$ satisfies {\bf (H$^{\alpha'}_{K'}$)} with $\alpha'=\alpha\wedge(\delta/2)$
and $K'=\kappa K$.
\item Let $\Theta^b$ be defined as above through $\Phi^b$. In the case of {\bf (3)}, we also have
\begin{align}
\|\Theta^b-\Theta^{b'}\|_{\mL^q_p(t_0,s_0)}\leq C\|b-b'\|_{\mL^q_p(t_0,s_0)},\label{KJ1}
\end{align}
where $C=C(K,\alpha,p,q,d,\delta,U)>0$.
\end{enumerate}
\el
\begin{proof}
Let $\delta:=\frac{1}{2}-\frac{d}{2p}-\frac{1}{q}>0$. By (\ref{ER1}), there is a $C_0=C_0(K,\alpha,p,q,d)>0$ such that for
all $[t_0,s_0]\subset[T,S]$,
$$
\|\nabla\u_t\|_{\sC^{\delta/2}}\leq C_0(s_0-t_0)^{\delta/3}\exp\left\{C_0(s_0-t_0)^{\delta q/3}\|b\|^q_{\mL^q_p(t_0,s_0)}\right\}\|b\|_{\mL^q_p(t_0,s_0)}.
$$
For given $U>0$, let us choose $\eps=\eps\big(\delta, q, C_0,U\big)>0$ small enough so that for all $s_0-t_0\leq\eps$ and  $\|b\|_{\mL^q_p(t_0,s_0)}\leq U$,
$$
\sup_{t\in[t_0,s_0]}\|\nabla\u_t\|_{\sC^{\delta/2}}\leq1/2.
$$
In particular, we have
$$
|\u_t(x)-\u_t(y)|\leq|x-y|/2,\  t\in[t_0,s_0],
$$
which then gives (\ref{EQ1}) by definition (\ref{Def1}).\\
\\
{\bf (1)} It is obvious from (\ref{EQ1}).
\\
\\
{\bf (2)} It follows from definition (\ref{Def1}) and the estimates (\ref{KH1}), (\ref{ER1}).
\\
\\
{\bf (3)} It follows from definition (\ref{Def1}) and the estimate (\ref{ER11}).
\\
\\
{\bf (4)} It follows by generalized It\^o's formula (see \cite{Kr} or \cite[Lemma 4.3]{Zh2} for more details).
\\
\\
{\bf (5)} By definition, we can write
\begin{align*}
\Theta^b_s(y)-\Theta^{b'}_s(y)&=[\nabla\Phi^{b}_s\cdot \sigma_s]\circ\Phi^{b,-1}_s(y)-[\nabla\Phi^{b}_s\cdot \sigma_s]\circ\Phi^{b',-1}_s(y)\\
&+[(\nabla\Phi^b_s-\nabla\Phi^{b'}_s)\cdot \sigma_s]\circ\Phi^{b',-1}_s(y)=:I_1(s,y)+I_2(s,y).
\end{align*}
For $I_1(s,y)$, by (\ref{ES2}) we have
$$
|I_1(s,y)|\leq C(\cM g_s(\Phi^{b,-1}_s(y))+\cM g_s(\Phi^{b',-1}_s(y)))|\Phi^{b,-1}_s(y)-\Phi^{b',-1}_s(y)|,
$$
where $g_s(x):=|\nabla[\nabla\Phi^{b}_s\cdot \sigma_s](x)|\in \mL^q_p(t_0,s_0)$ by {\bf (2)}, and $\cM g_s$ is the Hardy-Littlewood maximal function.
Noticing that
\begin{align*}
\sup_y|\Phi^{b,-1}_s(y)-\Phi^{b',-1}_s(y)|&=\sup_y|y-\Phi^{b',-1}_s\circ\Phi^{b}_s(y)|\leq\|\nabla\Phi^{b',-1}_s\|_\infty\|\Phi^{b'}_s-\Phi^{b}_s\|_\infty,
\end{align*}
by the change of variables, {\bf (3)} and \eqref{Es30}, we obtain
\begin{align*}
\|I_1\|_{\mL^q_p(t_0,s_0)}&\leq C\|\cM g_\cdot(\Phi^{b,-1}_\cdot)+\cM g_\cdot(\Phi^{b',-1}_\cdot)\|_{\mL^q_p(t_0,s_0)}\|\Phi^{b,-1}-\Phi^{b',-1}\|_\infty\\
&\leq C\|\cM g\|_{\mL^q_p(t_0,s_0)}\|b-b'\|_{\mL^q_p(t_0,s_0)}\leq C\|g\|_{\mL^q_p(t_0,s_0)}\|b-b'\|_{\mL^q_p(t_0,s_0)}.
\end{align*}
For $I_2(s,y)$, by the change of variables and {\bf (3)} again, we have
\begin{align*}
\|I_2\|_{\mL^q_p(t_0,s_0)}&\leq C\|\nabla\Phi^b_\cdot-\nabla\Phi^{b'}_\cdot\|_{\mL^q_p(t_0,s_0)}\leq C\|b-b'\|_{\mL^q_p(t_0,s_0)}.
\end{align*}
Combining the above calculations, we obtain (\ref{KJ1}).
\end{proof}
We are now in a position to give

\begin{proof}[Proof of Theorem \ref{Main}]
Let $\eps$ be as in Lemma \ref{Le51}. Fix $t_0\in[T,S)$ and $s_0\in(t_0,S)$ with
$$
s_0-t_0\leq\eps.
$$
Let us first prove the theorem on the time interval $[t_0,s_0]$.
By Lemma \ref{Le51} and Theorem \ref{Th0}, it is easy to see that {\bf (A)}, {\bf (B)} and {\bf (C)} hold.
Let us look at {\bf (D)}. By {\bf (d)} of Theorem \ref{Th0}, we have
\begin{align}
\nabla\mE f(Y_{t_0,s}(y))=\frac{1}{s-t_0}\mE\left(f(Y_{t_0,s}(y))\int^s_{t_0}\Theta^{-1}_r(Y_{t_0,r}(y))\nabla Y_{t_0,r}(y)\dif W_r\right).
\end{align}
Since $Y_{t_0,s}(y)=\Phi_s\circ X_{t_0,s}\circ \Phi^{-1}_{t_0}(y)$, by replacing $f$ with
$f\circ\Phi^{-1}_s$ and the change of variables $y\to\Phi_t(x)$,  we obtain (\ref{EU05}). As for {\bf (E)}, it follows by
{\bf (e)} of Theorem \ref{Th0} and (\ref{KJ1}).

Finally, let us consider the time interval $[t_1, s_1]$, where $t_1:=\frac{s_0+t_0}{2}$ and $s_1:=\frac{3s_0-t_0}{2}$.
By the uniqueness of solutions, we have for all $s\in[t_1,s_1]$,
$$
X_{t_0, s}(x)=X_{t_0, t_1}\circ X_{t_1,s}(x),
$$
where $X_{t_0,t_1}(\cdot)$ and $X_{t_1,s}(\cdot)$ are independent. Thus, we can patch up the solutions and conclude the proofs by Proposition \ref{Pr74}.
\end{proof}

\section{Proof of Theorem \ref{Th12}}

Given $p>d$, $\nu>0$ and $T\in[-1,0]$, let $b\in \mL^\infty_p(T,0)$ be divergence free, and let $X_{t,s}(x)$ solve
\begin{align}
X_{t,s}(x)=x+\int^s_t b_r(X_{t,r}(x))\dif r+\sqrt{2\nu}(W_s-W_t),\ T\leq t\leq s\leq 0.\label{SDE1}
\end{align}
\bl
For any $f\in L^1(\mR^d)$, we have
\begin{align}
\mE\int_{\mR^d}f(X_{t,s}(x))\dif x=\int_{\mR^d}f(x)\dif x.\label{EW2}
\end{align}
\el
\begin{proof}
By a density and monotonic class argument, it suffices to prove it for $f\in C^\infty_0(\mR^d)$.
Let $b^n_t(x)=\varrho_n*b_t(x)$, where $\rho_n$ is a mollifier. Then $\|\nabla b^n\|_\infty<\infty$ and $\div b^n_t=0$.
Since
$$
\det(\nabla X^n_{t,s}(x))=\exp\left\{\int^s_t\div b^n_r(X^n_{t,r}(x))\dif r\right\}=1,
$$
by the change of variables, one has
\begin{align}
\int_{\mR^d}f(X^n_{t,s}(x))\dif x=\int_{\mR^d}f(x)\det(\nabla X^{n,-1}_{t,s}(x))\dif x=\int_{\mR^d}f(x)\dif x,\label{BB1}
\end{align}
where $x\mapsto X^{n,-1}_{t,s}(x)$ is the inverse of $x\mapsto X^n_{t,x}(x)$.
On the other hand, by (\ref{EU55}) we have
$$
\lim_{n\to\infty}\mE\left(\sup_{s\in[t,0]}|X^n_{t,s}(x)-X_{t,s}(x)|^2\right)=0.
$$
By taking limits for both sides of \eqref{BB1}, we obtain \eqref{EW2}.
\end{proof}

Let $\bP=\mI-\nabla(-\Delta)^{-1}\div$ be Leray's projection onto the space of divergence free vector fields.
It is well-known that the singular integral operator  $\bP$ is bounded from $L^p$ to $L^p$ (cf. \cite[Theorem 3, p.96]{St}).
We also need the following result (cf. \cite{Co-Iy} and \cite{Zh1}).
\bl
Recall the definition of $\sV^0_{\infty-}$ in Section 2. Let $\varphi\in \mW^{1}_p(\mR^d;\mR^d)$ for some $p>1$. We have the following conclusions:
\begin{enumerate}[(i)]
\item For any $X\in L^\infty_x(L^{\infty-}_\omega)\cap\sV_{\infty-}$ and $Y\in\sV^0_{\infty-}$, we have
\begin{align}
\bP\mE[\nabla^{\mathrm{t}} X\cdot\varphi(Y)]=-\bP\mE[\nabla^{\mathrm{t}} Y\cdot\nabla^{\mathrm{t}}\varphi(Y)\cdot X].\label{FG1}
\end{align}
\item For any $X\in\sV^0_{\infty-}$, we have
\begin{align}
\nabla\bP\mE[\nabla^{\mathrm{t}} X\cdot\varphi(X)]=\bP\mE[\nabla^{\mathrm{t}} X\cdot(\nabla^{\mathrm{t}}\varphi-\nabla\varphi)(X)\cdot\nabla X].\label{FG2}
\end{align}
\end{enumerate}
\el
\begin{proof}
Let $X_n, Y_n,\varphi_n$ be the mollifying approximations of $X,Y,\varphi$  defined as in (\ref{R1}).
\\
(i) Notice that
\begin{align*}
&\bP\mE[\nabla^{\mathrm{t}} X_n\cdot\varphi_n(Y_m)]+\bP\mE[\nabla^{\mathrm{t}} Y_m\cdot\nabla^{\mathrm{t}}\varphi_n(Y_m)\cdot X_n]
=\bP\nabla\mE[X_n\cdot\varphi_n(Y_m)]=0.
\end{align*}
By (\ref{R2}), the dominated convergence theorem and H\"older's inequality, it is easy to see that for each $n\in\mN$,
$$
\mE[\nabla^{\mathrm{t}} X_n\cdot\varphi_n(Y_m)]\to
\mE[\nabla^{\mathrm{t}} X_n\cdot\varphi_n(Y)]\mbox{ in $L^p$ as $m\to\infty$},
$$
and
$$
\mE[\nabla^{\mathrm{t}} Y_m\cdot\nabla^{\mathrm{t}}\varphi_n(Y_m)\cdot X_n]\to
\mE[\nabla^{\mathrm{t}} Y\cdot\nabla^{\mathrm{t}}\varphi_n(Y)\cdot X_n]\mbox{ in $L^p$ as $m\to\infty$}.
$$
Hence,
$$
\bP\mE[\nabla^{\mathrm{t}} X_n\cdot\varphi_n(Y)]=-\bP\mE[\nabla^{\mathrm{t}} Y\cdot\nabla^{\mathrm{t}}\varphi_n(Y)\cdot X_n].
$$
By letting $n\to\infty$, we obtain (\ref{FG1}).
\\
\\
(ii) As above calculations, we have
$$
\nabla\bP\mE[\nabla^{\mathrm{t}} X_m\cdot\varphi_n(X_m)]
=\bP\mE[\nabla^{\mathrm{t}} X_m\cdot(\nabla^{\mathrm{t}}\varphi_n-\nabla\varphi_n)(X_m)\cdot\nabla X_m].
$$
By H\"older's inequality, we have
$$
\sup_{n,m}\|\nabla\bP\mE[\nabla^{\mathrm{t}} X_m\cdot\varphi_n(X_m)]\|_p<\infty.
$$
Firstly letting $m\to\infty$ and then $n\to\infty$, we find that
$$
\mE[\nabla^{\mathrm{t}} X_m\cdot(\nabla^{\mathrm{t}}\varphi_n-\nabla\varphi_n)(X_m)\cdot\nabla X_m]
\to\mE[\nabla^{\mathrm{t}} X\cdot(\nabla^{\mathrm{t}}\varphi-\nabla\varphi)(X)\cdot\nabla X]\mbox{ in $L^p$},
$$
and
$$
\mE[\nabla^{\mathrm{t}} X_m\cdot\varphi_n(X_m)]\to \mE[\nabla^{\mathrm{t}} X\cdot\varphi(X)]\mbox{ in $L^p$}.
$$
Combining the above calculations, we obtain (\ref{FG2}).
\end{proof}

Below we fix
$$
p>d \mbox{ and } q>(2p)/(p-d),
$$
and for given $\varphi\in L^p(\mR^d;\mR^d)$, define
$$
\mT(b)_t(x):=u_t(x):=\bP\mE[\nabla^{\mathrm{t}}X_{t,0}\cdot\varphi(X_{t,0})](x).
$$
\bl\label{Le62}
For any given $\varphi\in L^p(\mR^d)$, there exist a constant $C_0=C_0(d,p,q,\nu)>0$ and a time $T_0=T_0(C_0,\|\varphi\|_p)<0$
such that if $\|b\|_{\mL^\infty_p(T_0,0)}\leq 2C_0\|\varphi\|_p$ and $\div b=0$, then
$$
\|\mT(b)_t\|_{p}\leq 2C_0\|\varphi\|_p, \  t\in[T_0,0].
$$
\el
\begin{proof}
Let $\|\cdot\|_{L^p_{x,\omega}}$ be the norm in $L^p(\mR^d\times\Omega; \dif x\times P)$.
By definition and \eqref{EW2}, we have
\begin{align*}
\|\mT(b)_t\|_{p}&\leq C_{d,p}\|\mE[\nabla^{\mathrm{t}}X_{t,0}\cdot\varphi(X_{t,0})]\|_{p}\\
&\leq C_{d,p}\mathrm{ess.}\sup_{x\in\mR^d}\|\nabla^{\mathrm{t}}X_{t,0}(x)\|_{L^2_\omega}\|\varphi(X_{t,0})\|_{L^p_{x,\omega}}\\
&=C_{d,p}\mathrm{ess.}\sup_{x\in\mR^d}\|\nabla^{\mathrm{t}}X_{t,0}(x)\|_{L^2_\omega}\|\varphi\|_{L^p_x}\\
&\leq C(d, q, p, \nu, \|b\|_{\mL^q_p(t,0)})\|\varphi\|_{p},
\end{align*}
where the first inequality is due to the boundedness of {\bf P} in $L^p$, and the last inequality is due to {\bf (B)} of Theorem  \ref{Main}.
Since the constant $C$ is increasing with respect to $\|b\|_{\mL^q_p(t,0)}$ and goes to some $C_0=C_0(d,p,q,\nu)$
as $\|b\|_{\mL^q_p(t,0)}\to 0$, and also noticing that
$$
\|b\|_{\mL^q_p(t,0)}\leq \|b\|_{\mL^\infty_p(t,0)}|t|^{1/q}\leq 2C_0|t|^{1/q}\|\varphi\|_p,
$$
one can choose $T_0<0$ close to zero so that
$$
C(d, q, p, \nu, 2C_0|T_0|^{1/q}\|\varphi\|_p)\leq 2C_0.
$$
The proof is complete.
\end{proof}
\bl\label{Le72}
For given $\varphi\in \mW^1_p(\mR^d;\mR^d)$, let $C_0$ and $T_0$ be as in Lemma \ref{Le62}
and $U:=2C_0\|\varphi\|_{\mW^1_p}$, there exists a time $T_1=T_1(d,\nu,p,q,U)\in[T_0,0)$
such that for all $b,b'\in\mL^\infty_p(T_1,0)$ with
$$
\|b\|_{\mL^\infty_p(T_1,0)},\|b'\|_{\mL^\infty_p(T_1,0)}\leq U,\ \ \div b=\div b'=0,
$$
it holds that for all $t\in[T_1,0]$,
$$
\|\mT(b)_t-\mT(b')_t\|_{p}\leq \tfrac{1}{2}\|b-b'\|_{\mL^\infty_p(T_1,0)}.
$$
\el
\begin{proof}
Let $X^b_{t,0}$ be the solution of SDE \eqref{SDE1} with drift $b$. By definition, we have
\begin{align*}
\|\mT(b)_t-\mT(b')_t\|_{p}&\leq
\|\bP\mE(\nabla^{\mathrm{t}}X^b_{t,0}\cdot\varphi(X^b_{t,0}))-\bP\mE(\nabla^{\mathrm{t}}X^{b'}_{t,0}\cdot\varphi(X^{b'}_{t,0}))\|_{p}\\
&\leq \|\bP\mE(\nabla^{\mathrm{t}}X^{b'}_{t,0}\cdot(\varphi(X^b_{t,0})-\varphi(X^{b'}_{t,0})))\|_{p}\\
&+\|\bP\mE(\nabla^{\mathrm{t}}(X^b_{t,0}-X^{b'}_{t,0})\cdot\varphi(X^b_{t,0}))\|_p=:I_1+I_2.
\end{align*}
For $I_1$, by the boundedness of $\bP$ in $L^p$ and H\"older's inequality, we have
\begin{align}
\begin{split}
I_1&\leq C\|\mE(\nabla^{\mathrm{t}}X^{b'}_{t,0}\cdot(\varphi(X^b_{t,0})-\varphi(X^{b'}_{t,0})))\|_{p}\\
&\leq C\|\|\nabla^{\mathrm{t}}X^{b'}_{t,0}\|_{L^{p_1}_\omega}\cdot\|\varphi(X^b_{t,0})-\varphi(X^{b'}_{t,0})\|_{L^{p_2}_\omega}\|_{p},
\end{split}\label{EW1}
\end{align}
where $\frac{1}{p_1}+\frac{1}{p_2}=1$ with $p_2\in(1,\frac{2p}{p+2})$. By (\ref{ES2}) and {\bf (E)} of Theorem \ref{Main}, we have
\begin{align*}
\mE|\varphi(X^b_{t,0})-\varphi(X^{b'}_{t,0})|^{p_2}
&\leq C\mE\Big((\cM|\nabla\varphi|(X^b_{t,0})+\cM|\nabla\varphi|(X^{b'}_{t,0}))^{p_2}|X^b_{t,0}-X^{b'}_{t,0}|^{p_2}\Big)\\
&\leq C\Big(\mE(\cM|\nabla\varphi|(X^b_{t,0})+\cM|\nabla\varphi|(X^{b'}_{t,0}))^{\frac{2p_2}{2-p_2}}\Big)^{1-\frac{p_2}{2}}\Big(\mE|X^b_{t,0}-X^{b'}_{t,0}|^2\Big)^{\frac{p_2}{2}}\\
&\leq C\Big(\mE(\cM|\nabla\varphi|(X^b_{t,0})+\cM|\nabla\varphi|(X^{b'}_{t,0}))^{\frac{2p_2}{2-p_2}}\Big)^{1-\frac{p_2}{2}}\|b-b'\|^{p_2}_{\mL^q_p(t,0)}.
\end{align*}
Substituting this into (\ref{EW1}), and by {\bf (B)} of Theorem \ref{Main} and (\ref{EW2}), we obtain
\begin{align}
\begin{split}
I_1&\leq C\left(\int_{\mR^d}\mE(\cM|\nabla\varphi|(X^b_{t,0})+\cM|\nabla\varphi|(X^{b'}_{t,0}))^p\dif x\right)^{\frac{1}{p}}\|b-b'\|_{\mL^q_p(t,0)}\\
&\leq C\|\cM|\nabla\varphi|\|_p\|b-b'\|_{\mL^q_p(t,0)}\leq C\|\nabla\varphi\|_p|t|^{\frac{1}{q}}\|b-b'\|_{\mL^\infty_p(t,0)}.
\end{split}\label{FG7}
\end{align}
As for $I_2$, letting $p'=\frac{2p}{p-2}$, by (\ref{FG1}), H\"older's inequality, (\ref{EW2})  and (\ref{EU44}), we have
\begin{align*}
I_2&=\|\bP\mE(\nabla^{\mathrm{t}}X^b_{t,0}\cdot\nabla^{\mathrm{t}}\varphi(X^b_{t,0})\cdot(X^b_{t,0}-X^{b'}_{t,0}))\|_p\\
&\leq C\|\|X^b_{t,0}-X^{b'}_{t,0}\|_{L^2_\omega}\cdot\|\nabla\varphi(X^b_{t,0})\|_{L^p_\omega}\cdot\|\nabla X^b_{t,0}\|_{L^{p'}_\omega}\|_p\\
&\leq C\|b-b'\|_{\mL^q_p(t,0)}\|\nabla\varphi(X^b_{t,0})\|_{L^p(\mR^d\times\Omega)}\cdot\|\nabla X^b_{t,0}\|_{L^\infty_xL^{p'}_\omega}\\
&\leq C\|\nabla\varphi\|_p|t|^{\frac{1}{q}}\|b-b'\|_{\mL^\infty_p(t,0)},
\end{align*}
which, together with (\ref{FG7}), and letting $T_1\in[T_0,0)$ be small enough, yields the desired estimate.
\end{proof}

We are now in a position to give
\begin{proof}[Proof of Theorem \ref{Th12}]
By Lemmas \ref{Le62} and \ref{Le72}, the nonlinear operator $\mT$ is a contraction operator in the ball of $\mL^\infty_p(T_1,0)$
with radius $U=2C_0\|\varphi\|_{\mW^1_p}$. Therefore, by Banach's fixed point theorem, there is a unique point $u\in \mL^\infty_p(T_1,0)$ such that for
each $t\in[T_1,0]$,
$$
\mT(u)_t=u_t.
$$
On the other hand, by (\ref{FG2}), H\"older's inequality and (\ref{EU44}), (\ref{EW2}), we also have
$$
\|\nabla\mT(u)_t\|_p\leq C\|\mE[|\nabla X_{t,0}|^2\cdot|\nabla^{\mathrm{t}}\varphi-\nabla\varphi|(X_{t,0})]\|_p<+\infty.
$$
The proof is complete.
\end{proof}

{\bf Acknowledgements:}

Deep thanks go to the referee for his/her very carefully reading the manuscript and useful suggestions. Special thanks also go to Professor Renming Song for 
improving the writing.
This work is supported by NNSFs of China (Nos. 11271294, 11325105).

\end{document}